\documentclass{amsart}
\usepackage{amssymb,amsmath, amsthm,latexsym}
\usepackage{graphics}
\usepackage{amscd}
\usepackage{graphics}
\usepackage{here}
\newcommand{\cal}[1]{\mathcal{#1}}
\theoremstyle{plain}
\newtheorem{theorem}{Theorem}
\newtheorem{lemma}{Lemma}[section]
\newtheorem{theo}[lemma]{Theorem}
\newtheorem{proposition}[lemma]{Proposition}
\newtheorem{corollary}[lemma]{Corollary}
\parskip=\bigskipamount

\let\egthree=\phi
\let\phi=\varphi
\let\varphi=\egthree

\begin{document}
\title{Isometry groups of proper hyperbolic spaces}
\author{Ursula Hamenst\"adt}
\thanks 
{Partially supported by Sonderforschungsbereich 611.}

\date{February 23, 2008}

\begin{abstract}
Let $X$ be a proper hyperbolic geodesic metric space
and let $G$
be a closed subgroup of the isometry group ${\rm Iso}(X)$ of $X$.
We show that if $G$ is not elementary then 
for every $p\in (1,\infty)$ the second continuous
bounded cohomology group $H_{cb}^2(G,L^p(G))$ does not vanish.
As an application, we derive some
structure results for closed subgroups of ${\rm Iso}(X)$.
\end{abstract}

\maketitle

\section{Introduction}

A geodesic metric space $X$ is called
\emph{$\delta$-hyperbolic} for some
$\delta >0$ if it satisfies the $\delta$-thin triangle condition:
For every geodesic triangle in $X$ with sides $a,b,c$ the side $a$
is contained in the $\delta$-neighborhood of $b\cup c$.
If $X$ is \emph{proper} (i.e. if closed balls in $X$ of finite
radius are compact) then it
can be compactified by adding the
\emph{Gromov boundary} $\partial X$.
Moreover, the
isometry group ${\rm Iso}(X)$ of $X$, equipped with the compact
open topology, is a locally compact $\sigma$-compact
topological group which acts as a group of homeomorphisms
on $\partial X$.
The \emph{limit set}
$\Lambda$ of a subgroup $G$ of ${\rm Iso}(X)$ is
the set of accumulation points in $\partial X$ of an orbit
of the action of $G$ on $X$. This limit
set is a closed $G$-invariant subset of $\partial X$. 
The group $G$ is called
\emph{elementary} if its limit set consists of at most two
points.

A \emph{compact extension} of a
topological group $H$ is a topological group $G$ which contains a
compact normal subgroup $K$ such that $H=G/K$ as topological
groups. Extending earlier results of Monod-Shalom \cite{MS04} and of
Mineyev-Monod-Shalom \cite{MMS04}, 
we show.

\begin{theorem} Let $X$ be a proper hyperbolic geodesic
metric space and let $G< {\rm Iso}(X)$ be a
closed subgroup.
Then one of the following three
possibilities holds.
\begin{enumerate}
\item $G$ is elementary.
\item Up to passing to an open subgroup of finite
index, $G$ is a compact extension of a simple Lie group of rank one.
\item $G$ is a compact extension of a totally disconnected group.
\end{enumerate}
\end{theorem}

The proof of the above result uses continuous second bounded
cohomology for locally compact topological groups $G$
with coefficients in a 
\emph{Banach module} for $G$. Such a Banach module 
is a separable Banach space $E$ together with a
continuous
homomorphism of $G$ into the group of linear
isometries of $E$. For every such Banach module $E$
for $G$ and every $i\geq 1$, the
group $G$ naturally acts on the vector space
$C_b(G^i,E)$ of continuous bounded maps $G^i\to E$.
If we denote by
$C_b(G^i,E)^G\subset C_b(G^i,E)$
the linear subspace of all $G$-invariant such maps, then the
\emph{second continuous bounded
cohomology group} $H_{cb}^2(G,E)$ of $G$
with coefficients $E$
is defined as the second cohomology group of the complex
\begin{equation} 0\to
C_b(G,E)^G \xrightarrow{d} C_b
(G^2,E)^G\xrightarrow{d} \dots \notag
\end{equation} with the usual
homogeneous coboundary operator $d$ (see \cite{M}).
The trivial representation of $G$ on $\mathbb{R}$ defines
the real second bounded cohomology group $H_{cb}^2(G,\mathbb{R})$.
Let $H_c(G,\mathbb{R})$ be the ordinary 
\emph{continuous} real second cohomology group 
of $G$. Then there is a natural
homomorphism $H_{cb}(G,\mathbb{R})\to H_c(G,\mathbb{R})$
which in general is neither injective nor surjective.

A closed subgroup $G$ of ${\rm Iso}(X)$ is a locally compact and
$\sigma$-compact topological group 
and hence it admits a left invariant locally
finite Haar measure $\mu$. In particular, for every $p\in (1,\infty)$
the separable Banach space $L^p(G,\mu)$ of
functions on $G$ which are $p$-integrable with respect to $\mu$
is a Banach module for $G$ 
with respect to the isometric action of $G$ by left translation.
Extending earlier results of 
Monod-Shalom \cite{MS04}, of 
Mineyev-Monod-Shalom \cite{MMS04} and of Fujiwara \cite{F98}, 
we obtain the following non-vanishing result
for second bounded cohomology.

\begin{theorem} 
Let $G$ be a closed non-elementary subgroup of the isometry
group of a proper hyperbolic geodesic metric space $X$ with
limit set $\Lambda\subset \partial X$. Then we have
$H_{cb}^2(G,L^p(G,\mu))\not=\{0\}$ for every $p\in (1,\infty)$.
Moreover, one of the following two possibilities holds.
\begin{enumerate}
\item 
$G$ does not act transitively on the complement of the diagonal
in $\Lambda\times \Lambda$ and the kernel of the natural homomorphism
$H_{cb}^2(G,\mathbb{R})\to H_c^2(G,\mathbb{R})$ is infinite dimensional.
\item $G$ acts transitively on the complement of the 
diagonal in $\Lambda\times \Lambda$ and the
kernel of the natural homomorphism
$H_{cb}^2(G,\mathbb{R})\to H_c^2(G,\mathbb{R})$ is trivial.
\end{enumerate}
\end{theorem}

Theorem 2 implies an extension of all the results of
\cite{MS04,MMS04} to arbitrary proper hyperbolic geodesic metric
spaces without any additional assumptions. 

We also investigate cocycles with values in ${\rm Iso}(X)$.
Namely,
let $S$ be a standard Borel space and let $\nu$ be a Borel
probability measure on $S$. Let $\Gamma$ be a countable group
which admits a measure preserving action on $(S,\nu)$. 
The action is called \emph{mildly mixing} if there are no
nontrivial recurrent sets, i.e. if for every
Borel subset
$A$ of $S$ with $\nu(A)\in (0,1)$ and every sequence
$g_i\to \infty$ in $\Gamma$ we have
$\lim\inf_{i\to \infty}\nu(A\Delta g_iA)\not= 0$. An
${\rm Iso}(X)$-valued \emph{cocycle} for the action is a
measurable map $\alpha:\Gamma\times S\to {\rm Iso}(X)$ such that
\begin{equation}
\alpha(gh,x)=\alpha(g,hx)\alpha(h,x) \notag
\end{equation} 
for all $g,h\in \Gamma$ and for 
$\nu$-almost every $x\in S$. The cocycle $\alpha$ is
\emph{cohomologous} to a cocycle $\beta:\Gamma\times S\to {\rm
Iso}(X)$ if there is a measurable map $\phi:S\to {\rm Iso}(X)$
such that \begin{equation}
\phi(gx)\alpha(g,x)=\beta(g,x)\phi(x) \notag
\end{equation} 
for all $g\in G$, $\nu$-almost every $x\in S$. 

A lattice $\Gamma$ in a product $G_1\times G_2$ of two
locally compact $\sigma$-compact 
and non-compact topological groups is
called \emph{irreducible} if the projection of $\Gamma$ into
each factor is dense.
Extending
earlier results of Monod and Shalom \cite{MS04} 
we show.

\begin{theorem} Let $G$ be a semi-simple Lie group with
finite center, no compact factors and of rank at least $2$.
Let $\Gamma<G$ be an irreducible lattice which admits a mildly 
mixing measure
preserving action on a standard probability space $(S,\nu)$. Let $X$ be a
proper hyperbolic geodesic metric space and let
$\alpha:\Gamma\times S\to {\rm Iso}(X)$ be a cocycle. Then one of
the following two possibilities holds.
\begin{enumerate}
\item $\alpha$ is cohomologous to a cocycle into
an elementary subgroup of ${\rm Iso}(X)$.
\item $\alpha$ is cohomologous to a cocycle
into a subgroup $H$ of ${\rm Iso}(X)$
which is a compact extension of a simple
Lie group $L$ of rank one,
and there is a continuous surjective homomorphism $G\to L$.
\end{enumerate}
\end{theorem}

A proper hyperbolic geodesic metric space $X$ is said to be of
\emph{bounded growth} if there is a constant $b>0$ such that for
every $R\geq 1$ an open metric ball in $X$ of radius $R$ contains
at most $b e^{bR}$ disjoint open metric balls of radius $1$. Every
hyperbolic metric graph of bounded valence and every finite
dimensional simply connected Riemannian manifold of bounded 
negative curvature
has this property. We
show.

\begin{theorem} Let $\Gamma$ be a finitely generated group
which admits a properly discontinuous isometric action on a
hyperbolic geodesic metric space $X$ of bounded growth. If
$H_{cb}^2(\Gamma,\mathbb{R})$ or
$H_{cb}^2(\Gamma,\ell^2(\Gamma))$ is finite dimensional then $\Gamma$
is virtually nilpotent.
\end{theorem}

The organization of this paper is as follows. In Section 2 we construct
for every non-elementary closed subgroup $G$ of the isometry
group
of a proper hyperbolic geodesic metric space $X$ with limit set
$\Lambda\subset \partial X$ and for every $p\in (1,\infty)$ 
a nontrivial
continuous bounded $L^p(G\times G,\mu\times \mu)$-valued
cocycle for the action of $G$ on $\Lambda$ where
as before, $\mu$ is a Haar measure on $G$. 
We use this in Section 3
to show that for every $p\in (1,\infty)$ the
second bounded cohomology group
$H_{cb}^2(G,L^p(G,\mu))$ does not vanish. 
In Section 4 we deduce
Theorem 1 from this non-vanishing result. 
The proof of Theorem 3 is contained in Section 5.

In Section 6 we construct for a closed non-elementary subgroup $G$ of
${\rm Iso}(X)$ with limit set $\Lambda\subset \partial X$ 
which does not act
transitively on the complement of the diagonal in 
$\Lambda\times \Lambda$ 
infinitely many quasi-morphisms whose linear spans are pairwise
inequivalent. This is then used to derive Theorem 2. The proof of 
Theorem 4 is contained in Section 7.

\section{Continuous bounded cocycles}

In this section let $(X,d)$ always be a
proper hyperbolic geodesic metric space.
Recall that the \emph{Gromov boundary} $\partial X$ of $X$
is defined as follows. Fix a point $x\in X$.
For two points $y,z\in X$, define
the \emph{Gromov product} $(y,z)_{x}$
based at $x$ by
\begin{equation}
(y,z)_{x}=\frac{1}{2}\bigl(d(y,x)+d(z,x)-d(y,z)\bigr). \notag
\end{equation}

Let ${\cal A}\subset \prod_{i\geq 0}X$ be the set of
all sequences $(y_i)\subset X$ such that
$(y_i,y_j)_x\to \infty$ $(i,j\to \infty)$.
Call two sequences $(y_i),(z_j)\in {\cal A}$ \emph{equivalent} if
$(y_i,z_i)_x\to \infty$ $(i\to \infty)$. By hyperbolicity of $X$,
this notion of equivalence defines an equivalence relation on
${\cal A}$ \cite{BH}.
The boundary $\partial X$ of $X$ is the set of equivalence classes
of this relation.

The Gromov product $(\,,\,)_{x}$ for pairs of points in $X$
can be extended to a product on
$\partial X$ by defining \begin{equation} (\xi,\eta)_{x}= \sup
\liminf_{i,j\to\infty}(y_i,z_j)_{x}\notag
\end{equation} where the supremum is taken over all sequences
$(y_i),(z_j)\subset X$
whose equivalence classes define the
points $\xi,\eta\in
\partial X$. 

For a number $\chi>0$
only depending on the hyperbolicity constant of $X$, there is a
distance $\delta_{x}$ on $\partial X$
with the property that the distance $\delta_x(\xi,\eta)$ between two
points $\xi,\eta\in \partial X$ is comparable to
$e^{-\chi(\xi,\eta)_{x}}$ (see 7.3 of \cite{GH}). More precisely,
there is a constant $\theta>0$ not depending on $x$ such that
\begin{equation}\label{dist}
e^{-\chi\theta} e^{-\chi(\xi,\eta)_{x}}\leq \delta_x(\xi,\eta)\leq
e^{-\chi(\xi,\eta)_{x}}\end{equation} for all $\xi,\eta\in
\partial X$. 
In particular, the diameter of $\delta_x$ is bounded
from above independent of $x$.
The distances $\delta_x$ $(x\in X)$ on $\partial X$  
are invariant under the natural action of
the isometry group ${\rm Iso}(X)$ 
of $X$ on $\partial X\times \partial X\times X$,
and they satisfy
\begin{equation}\label{distequi}
e^{-\chi d(x,y)}\delta_x\leq \delta_y\leq e^{\chi d(x,y)}\delta_x
\text{ for all } x,y\in X.
\end{equation}
As a consequence, the topology on $\partial X$ defined by the metric
$\delta_x$ does not depend on $x$. 

There is a natural topology on $X\cup \partial X$ which restricts
to the given topology on $X$ and to the topology on $\partial X$
induced by any one of the metrics $\delta_x$.
With respect to this topology, a
sequence $(y_i)\subset X$ converges to $\xi\in
\partial X$ if and only if the sequence is contained in ${\cal A}$
(i.e. we have $(y_i,y_j)_{x}\to \infty$) and if moreover the
equivalence class of $(y_i)$ equals $\xi$. Since $X$ is proper
by assumption, the space
$X\cup \partial X$ is compact and metrizable.
Every isometry of $X$ acts
naturally on $X\cup
\partial X$ as a homeomorphism (see \cite{BH} for all this).

We need the following simple observation.

\begin{lemma}\label{productmetric}
$\partial X\times X$ admits a natural ${\rm Iso}(X)$-invariant
distance function $\tilde d$ inducing the product topology.

\end{lemma}
\begin{proof}
Let $\delta_x$ $(x\in X)$
be a family of metrics on $\partial X$ which is
invariant under the action of ${\rm Iso}(X)$ on
$\partial X\times \partial X\times X$ and 
such that for some $\theta>0,\chi \in (0,1)$ the inequalities
(\ref{dist}) and (\ref{distequi}) above are valid. 

For $(\xi,x),(\eta,y)\in \partial X\times X$ write 
\begin{equation}\label{tilded}
\tilde
d_0((\xi,x),(\eta,y))=d(x,y)+\frac{1}{2}(\delta_x(\xi,\eta)+
\delta_y(\xi,\eta)).\notag
\end{equation}
Note that $\tilde d_0$ is a symmetric function
on $(\partial X\times X)\times (\partial X\times X)$
which is invariant
under the diagonal action of ${\rm Iso}(X)$. 

Define a
\emph{chain} between $(\xi,x),(\eta,y)$ to be a finite sequence
$\{(\xi_i,x_i)\mid 0\leq i\leq m\}\subset \partial X\times X$ with
$(\xi_0,x_0)=(\xi,x)$ and $(\xi_m,x_m)=(\eta,y)$. The \emph{length} of the
chain is defined to be $\sum_i\tilde d_0((\xi_i,x_i),
(\xi_{i+1},x_{i+1}))$. Let $\tilde d((\xi,x),(\eta,y))$ be the infimum of
the lengths of any chains connecting $(\xi,x)$ to $(\eta,y)$. 
Then $\tilde d$ is a symmetric ${\rm Iso}(X)$-invariant
function on $(\partial X\times X)\times (\partial X\times X)$ 
which is bounded from above by 
$\tilde d_0.$  By construction, $\tilde d$ satisfies
the triangle inequality.
Moreover, by the triangle inequality for $d$ and 
the definition, the function
$\tilde d$ satisfies
\begin{equation}
\tilde d((\xi,x),(\eta,y))\geq d(x,y)\,\forall\xi,\eta\in \partial X,
\forall x,y\in X.\notag\end{equation}

We claim that there is a constant $c>0$ such that
\begin{equation}\label{tildedelta}
c\delta_x(\xi,\eta)\leq \tilde d((\xi,x),(\eta,x))\leq \delta_x(\xi,\eta)
\,\forall x\in X,\forall \xi,\eta\in \partial X.
\end{equation}
Namely, the upper estimate for $\tilde d((\xi,x),(\eta,x))$ is
immediate from the definitions. To show the lower estimate, note that
the diameter of $\partial X$ equipped with any one of the metrics
$\delta_x$ is bounded from above by a universal constant $D>0$.
Let $x\in X$, let $\xi,\eta\in \partial X$ and let 
$\{(\xi_i,x_i)\mid 0\leq i\leq m\}$ be a chain
connecting $(\xi,x)=(\xi_0,x_0)$ to $(\eta,x)=(\xi_m,x_m)$.
If there is some $i<m$ such that $d(x,x_i)\geq D$ then
by the triangle inequality for $d$,  
the length of the chain is bigger than $2D$ and hence it is
bigger than
$\tilde d((\xi,x),(\eta,x))+D$.
On the other hand, if for every $i\leq m$ we have $d(x,x_i)\leq D$
then inequality (\ref{distequi}) together with the triangle inequality
for $\delta_x$ 
shows that the length of the chain is not smaller than
$e^{-\chi D}\delta_x(\xi,\eta)$. 
Inequality (\ref{tildedelta}) follows..

As a consequence, we have $\tilde d((\xi,x),(\eta,y))=0$ only
if $\xi=\eta,x=y$. This shows that 
$\tilde d$ is a distance on
$\partial X\times X$ which moreover induces the product 
topology.
\end{proof}

For every proper metric space $X$, the isometry group ${\rm
Iso}(X)$ of $X$ can be equipped with a natural locally compact
$\sigma$-compact metrizable topology, the so-called \emph{compact
open topology}. With respect to this topology, a sequence
$(g_i)\subset {\rm Iso}(X)$ converges to some isometry $g$ if and
only if $g_i\to g$ uniformly on compact subsets of $X$.
In this topology, a
closed subset $A\subset {\rm Iso}(X)$ is compact if and only if
there is a compact subset $K$ of $X$ such that $gK\cap
K\not=\emptyset$ for every $g\in A$. In particular, the action of
${\rm Iso}(X)$ on $X$ is proper.
In the
sequel we always equip subgroups of ${\rm Iso}(X)$
with the compact open topology.

Let again $X$ be a proper hyperbolic geodesic metric space
and let $G<{\rm Iso}(X)$ be a subgroup of the isometry group of $X$.
The \emph{limit set} $\Lambda$ of $G$ is the set
of accumulation points in $\partial X$ of one (and hence every)
orbit of the action of $G$ on $X$. If the closure of 
$G$ is non-compact then its
limit set is a compact non-empty $G$-invariant subset of $\partial
X$. The group $G$ is called \emph{elementary} if its limit set
consists of at most two points. In particular, every compact
subgroup of ${\rm Iso}(X)$ is elementary. If $G$ is non-elementary
then its limit set $\Lambda$ is uncountable without isolated
points. 

An element $g\in {\rm Iso}(X)$ is called \emph{hyperbolic}
if it generates an infinite cyclic subgroup $G<{\rm Iso}(X)$ whose
limit set $\Lambda$ consists of two points $a,b\in \partial X$
which are fixed points for $g$. Moreover, $a$ is 
an \emph{attracting fixed point},
$b$ is a \emph{repelling fixed point} and $g$ 
acts with \emph{north-south-dynamics} with respect to these fixed
points. This means that for any two neighborhoods $U$ of $a$,
$V$ of $b$ there is some $k>0$ such that $g^k(\partial X-V)\subset U$
and $g^{-k}(\partial X-U)\subset V$.
If $G<{\rm Iso}(X)$ is non-elementary then the
fixed points of hyperbolic elements in $G$ form a dense subset of
$\Lambda$ (8.26 of \cite{GH}).

Since $X$ is proper, any two points $\xi\not=\eta\in \partial X$
can be connected by a geodesic in $X$. We need the following
simple observation for which we did not find a precise
reference in the literature.

\begin{lemma}\label{fixedpoint}
Let $X$ be a proper $\kappa_0$-hyperbolic geodesic metric space
for some $\kappa_0>0$. Let 
$g$ be a hyperbolic isometry of $X$
and let $x_0$ be a point on a geodesic connecting 
the two fixed points $a\not=b\in \partial X$
for the action of $g$ on $\partial X$.
If $k>0$ is such that $d(x_0,g^kx_0)\geq 4\kappa_0$
then there is no isometry $h$ of $X$ with $hx_0=x_0,
hg^kx_0=g^kx_0$ and $h(a)=b,h(b)=a$.
\end{lemma}
\begin{proof}
Let $X$ be $\kappa_0$-hyperbolic for some $\kappa_0>0$,
let $g\in {\rm Iso}(X)$ be a hyperbolic isometry and let
$a\not=b\in\partial X$ be the attracting and the repelling
fixed point for the action of 
$g$ on $\partial X$, respectively. Let $\gamma$ be a 
geodesic in $X$ connecting $b$ to $a$. 
Such a geodesic may not be unique, but 
the Hausdorff distance
between $\gamma$ and any other 
such geodesic is bounded from above by $\kappa_0$.
Since $g$ maps the geodesic $\gamma$ connecting $b$ to $a$
to a geodesic $g\gamma$ connecting $gb=b$ to $ga=a$, 
for every $x\in \gamma$ the orbit of $x$ under 
the infinite cyclic subgroup of
${\rm Iso}(X)$ generated by $g$ 
is contained in the
$\kappa_0$-neighborhood of $\gamma$.

Let $x_0\in \gamma$ and let $k>0$ be such that 
$d(x_0,g^kx_0)\geq 4\kappa_0$. Assume that $\gamma$ is
parametrized in such a way that $\gamma(0)=x_0$.
Since $\gamma(t)$ converges as $t\to\infty$ to the
attracting fixed point $a$ for the action of $g$ on
$\partial X$, we have 
$d(g^kx_0,\gamma[3\kappa_0,\infty))\leq \kappa_0$. 
Assume to the contrary that there is an isometry $h$ of $X$ which 
fixes both $x_0$ and $g^kx_0$ and exchanges $a$ and $b$.
Then $h$ maps the geodesic ray $\gamma[0,\infty)$ connecting 
$x_0$ to $a$ to a
geodesic ray $h\gamma[0,\infty)$ 
connecting $x_0$ to $b$. Thus the   
Hausdorff distance between $h\gamma[0,\infty)$ and
$\gamma(-\infty,0]$ is at most $\kappa_0$.
As a consequence, the $\kappa_0$-neighborhood of the 
ray $h\gamma[3\kappa_0,\infty)$ does not
intersect the $\kappa_0$-neighborhood of $\gamma[3\kappa_0,\infty)$.
However, the $\kappa_0$-neighborhood
of $\gamma[3\kappa_0,\infty)$ contains $g^kx_0=hg^kx_0$,
and since $h$ is an isometry,
the point $hg^kx_0$ is contained in the
$\kappa_0$-neighborhood of $h\gamma[3\kappa_0,\infty)$ as
well. This is a contradiction and shows the lemma.
\end{proof}

The following theorem is the main technical result of this paper.
It gives a simplification and extension of some results of
Monod and Shalom \cite{MS04} and of Mineyev, Monod and Shalom
\cite{MMS04} avoiding  the difficult concept of homological
bicombing. For its formulation, let again $X$ be a proper hyperbolic
geodesic metric space and let $G$ be a \emph{closed} subgroup of the
isometry group of $X$ with limit set $\Lambda$. Note that
$G$ is a locally compact $\sigma$-compact topological group.
Assume that $G$ is non-elementary and let
$T\subset \Lambda^3$ be the space of triples of pairwise distinct
points in $\Lambda$. Then $T$ is an uncountable topological space without
isolated points. The group $G$ acts diagonally on $T$ as a group
of homeomorphisms. 

For a Banach-module $E$ for $G$  
define an \emph{$E$-valued continuous
bounded two-cocycle} for the action of $G$ on $\Lambda$ to be a
continuous \emph{bounded} $G$-equivariant map $\omega:T\to E$
which satisfies the following two properties.
\begin{enumerate}
\item 
For every permutation $\sigma$ of
the three variables, the \emph{anti-symmetry condition}
$\omega\circ
\sigma={\rm sgn}(\sigma)\omega$ holds.
\item For every quadruple $(x,y,z,w)$
of distinct points in $\Lambda$, the \emph{cocycle equality}
\begin{equation}\label{coc}
\omega(y,z,w)-\omega(x,z,w)+\omega(x,y,w)-\omega(x,y,z)=0
\end{equation}
is satisfied.
\end{enumerate}

Recall that every locally compact $\sigma$-compact topological group $G$ 
admits a left invariant locally finite Haar measure $\mu$. 
For $p\in (1,\infty)$ denote by $L^p(G\times G,\mu\times \mu)$ 
the Banach space of all functions on $G\times G$ which
are $p$-integrable with respect to the product measure $\mu\times\mu$.
The group $G$ acts continuously and isometrically on 
$L^p(G\times G,\mu\times \mu)$ by left translation.
We
have.

\begin{theo}\label{cocycle}
Let $X$ be a proper hyperbolic geodesic
metric space and let $G< {\rm Iso}(X)$ be a closed non-elementary
subgroup. For every $p\in (1,\infty)$ and every triple $(a,b,\xi)\in
T$ of pairwise distinct points in $\Lambda$ such that $(a,b)$ is
the pair of fixed points of a hyperbolic element of $G$ there is
an $L^p(G\times G,\mu\times \mu)$-valued continuous bounded
cocycle $\omega$ for the action of $G$ on $\Lambda$ with
$\omega(a,b,\xi)\not=0$.
\end{theo}
\begin{proof}
Let $G< {\rm Iso}(X)$ be a non-elementary closed
subgroup of ${\rm Iso}(X)$. 
We divide the proof of the theorem into five steps.

{\sl Step 1:}

Let $x_0\in X$ be an arbitrary point
and denote by $G_{x_0}$ the stabilizer of $x_0$ in $G$. Then
$G_{x_0}$ is a compact subgroup of $G$, and the quotient space
$G/G_{x_0}$ is $G$-equivariantly 
homeomorphic to the orbit $Gx_0\subset X$ of $x_0$.
Note that $Gx_0$ is a closed subset of $X$ and hence it is locally compact.

Since $G$ is non-elementary by assumption, the limit set
$\Lambda\subset \partial X$ of $G$ is an uncountable closed subset
of $\partial X$ without isolated points \cite{BH}. 
In particular, $\Lambda$ is compact.
The group $G$
acts on the locally compact space
$\Lambda\times X$ as a group of homeomorphisms.
Denote by
$\Delta$ the diagonal in $\partial X\times
\partial X$. The
${\rm Iso}(X)$-invariant metric $\tilde d$ on $\partial X\times X$ 
constructed in Lemma \ref{productmetric}
induces a
$G$-invariant metric on $\Lambda\times G/G_{x_0}$ and hence a
$G$-invariant symmetrized product metric $\hat d$ on 
\begin{equation}V=(\Lambda\times
\Lambda-\Delta)\times G/G_{x_0}\times G/G_{x_0} \notag
\end{equation} 
by defining  
\begin{align}\label{productmetric2}
\hat d((\xi,\eta,x,y),& (\xi^\prime,\eta^\prime,x^\prime,y^\prime))\\
=\frac{1}{2}\bigl(\tilde d((\xi,x),(\xi^\prime,x^\prime))+
\tilde d((\eta,y),(\eta^\prime,y^\prime))+ &
\tilde d((\eta,x),(\eta^\prime,x^\prime))+
\tilde d((\xi,y),(\xi^\prime,y^\prime))\bigr).\notag
\end{align}
The distance $\hat d$ is also
invariant under the involution
$\iota:(\xi,\eta,x,y)\to (\eta,\xi,x,y)$ exchanging the first
two factors. Moreover, $\hat d$ induces the product topology on $V$.

Since $V$ is locally
compact and $G<{\rm Iso}(X)$ is closed,
the space $W=G\backslash V$ admits a natural metric $d_0$ as
follows. Let 
\begin{equation}P:V\to W \notag
\end{equation}
be the canonical projection and define
\begin{equation}\label{d0}
d_0(x,y)=\inf\{\hat d(\tilde x,\tilde y)\mid P\tilde x=x, P\tilde
y=y\}. \end{equation}
The topology induced by this metric is the quotient
topology for the projection $P$. In particular, $W$ is a
locally compact metric space. Since the action of $G$ commutes
with the isometric involution $\iota$, the map  
$\iota$ descends to an isometric involution 
of the metric space $(W,d_0)$
which we denote again by $\iota$.  

Let $\kappa_0>0$ be a hyperbolicity constant for $X$. 
Let as before $\delta_z$ $(z\in X)$ be a family of distance
functions on $\partial X$ which satisfy the properties
(\ref{dist}) and (\ref{distequi}) in the beginning of this section.
There is a universal constant $c_0>0$ depending only
on $\kappa_0$ with the following
property. For any $\xi\not= \eta\in
\partial X$, for every geodesic $\gamma$ connecting $\eta$ to $\xi$ and
for every point $z\in X$ which is
contained in the $\kappa_0$-neighborhood
of $\gamma$ we have $\delta_z (a,b)\geq 2c_0$. Moreover, there is 
a constant $\kappa >6\kappa_0$ only depending on 
$\kappa_0$ (and not on $\xi,\eta$) such that 
the set $\{z\in X\mid
\delta_z(\xi,\eta)\geq c_0\}$ is contained in the $\kappa$-tubular
neighborhood of any geodesic $\gamma$ connecting $\eta$ to $\xi$.

Choose a smooth function $\phi_0:[0,\infty)\to [0,1]$ such that
$\phi_0(t)=0$ for all $t\leq c_0$ and $\phi_0(t)=1$ for all
$t\geq 2c_0$.
Define a continuous function
$\phi:V\to [0,1]$ by
\begin{equation}\label{phi}
\phi(\xi,\eta,gx_0,hx_0)=
\phi_0(\delta_{gx_0}(\xi,\eta))\phi_0(\delta_{hx_0}(\xi,\eta)). 
\end{equation}
Note that $\phi$ is invariant under the diagonal action of $G$ and
under the involution $\iota$.
In particular, the support $V_0$ of $\phi$ is invariant under $G$, 
and it projects to a $\iota$-invariant 
closed subset $W_0=G\backslash V_0$ of $W$.

Let $\kappa_1> 3\kappa$ and 
let $\psi:[0,\infty)\to [0,1]$ be a smooth function which
satisfies $\psi(t)=1$ for $t\leq \kappa_1$ and $\psi(t)=0$ for
$t\geq 2\kappa_1$. Define a continuous 
function $\zeta:V\to [0,1]$ by
\begin{equation}\label{zeta}
\zeta(\xi,\eta,ux_0,hx_0)= \phi(\xi,\eta,ux_0,hx_0)\psi(d(ux_0,hx_0)).
\end{equation}
Note that $\zeta$ is invariant under the diagonal action of $G$
and under the involution $\iota$. 
In particular, $\zeta$ projects to a continuous function $\zeta_0$ on 
$W=G\backslash V$. 

{\sl Step 2:}

In equation (\ref{d0})
in Step 1 above, we defined a distance $d_0$ on the space
$W=G\backslash V$.
With respect to this distance, 
the involution $\iota$ acts non-trivially
and isometrically. Choose a small closed metric 
ball $B$ in $W$
which is disjoint from its image under $\iota$. 
In Step 5 below we 
will construct explicitly 
such balls $B$, however 
for the moment, we simply assume that such a ball exists.

Let ${\cal H}$ be
the vector space of all H\"older continuous functions $f:W\to
\mathbb{R}$ supported in $B$. 
An example of such a function can be obtained as follows.

Let $z$ be an interior point of $B$ and let $r>0$ be sufficiently
small that the closed metric 
ball $B(z,r)$ of radius $r$ about $z$ is contained
in $B$. Choose a smooth function $\alpha:\mathbb{R}\to [0,1]$ such that
$\alpha(t)=1$ for $t\in [r/2,\infty)$ and $\alpha(t)=0$ for 
$t\in (-\infty,0]$ and 
define $f(y)=\alpha(r-d_0(y,z))$. Since the function $y\to r-d_0(y,z)$ 
is one-Lipschitz on $W$ and $\alpha$ is smooth, the function 
$f$ on $W$ is Lipschitz, does not vanish at $z$ 
and is supported in $B$.  

Since $B$ is disjoint from $\iota(B)$ by assumption and since
$\iota$ is an isometry, every function $f\in {\cal H}$ admits
a natural extension to a H\"older 
continuous function $f_0$ 
on $W$ supported in $B\cup \iota(B)$ whose restriction to
$B$ coincides with the restriction of $f$ and 
which satisfies $f_0(\iota z)=-f_0(z)$ for
all $z\in W$. 
The function $\hat f=f_0\circ P:V\to \mathbb{R}$
is invariant under the action of $G$, and it 
is \emph{anti-invariant} under the involution 
$\iota$ of $V$, i.e. it satisfies
$\hat f(\iota(v))=-\hat f(v)$ for all $v\in V$ (here as before,
$P:V\to W$ denotes the canonical projection).

Equip $\tilde V=(\Lambda\times \Lambda -\Delta)\times G\times G$
with the product topology. The group $G$
acts on $G\times G$ by left translation, and it acts 
diagonally on $\tilde V$.
Denote again by $\iota$ the involution of $\tilde V$ exchanging
the first two factors.
There is a natural continuous 
projection $\Pi:\tilde V\to V$ which is
equivariant with respect to the action of $G$ and with respect to
the action of the involution $\iota$ on $\tilde V$ and $V$.
The function $\hat f$ on $V$ lifts to
a $G$-invariant $\iota$-anti-invariant
continuous function $\tilde f=\hat f\circ \Pi$ on
$\tilde V$.
The lift $\tilde \zeta=\zeta\circ \Pi$ of
the function $\zeta$ defined in equation (\ref{zeta})
in Step 1
is $G$-invariant and $\iota$-invariant.

For $\xi\not= \eta\in \Lambda$ write $F(\xi,\eta)=\{(\xi,\eta,z)\mid
z\in G\times G\}$. The sets $F(\xi,\eta)$ define a $G$-invariant
foliation ${\cal F}$ of $\tilde V$. The leaf 
$F(\xi,\eta)$ of ${\cal F}$
can naturally be identified with $G\times G$. 
For all $\xi\not=\eta\in
\Lambda$ and every function $f\in {\cal H}$ we denote
by $f_{\xi,\eta}$ 
the restriction of the function $\tilde f$ to $F(\xi,\eta)$,
viewed as a continuous function on
$G\times G$. Similarly, define $\zeta_{\xi,\eta}$ to be the restriction
of the function $\tilde\zeta$ to $F(\xi,\eta)$.
For every $f\in {\cal H}$, all $\xi\not=\eta\in \Lambda$ 
and all $g\in G$
we then have $f_{g\xi,g\eta}\circ g=f_{\xi,\eta}=-f_{\eta,\xi}$, moreover
$\zeta_{g\xi,g\eta}\circ g=\zeta_{\xi,\eta}=\zeta_{\eta,\xi}$.

Recall that the limit set $\Lambda$
of $G$ is uncountable without
isolated points. For $f\in {\cal H}$ and for an ordered
triple $(\xi,\eta,\beta)$ of pairwise distinct points in $\Lambda$ define
\begin{equation}\label{omega}
\omega(\xi,\eta,\beta)= f_{\xi,\eta} \zeta_{\xi,\eta}+ 
f_{\eta,\beta}\zeta_{\eta,\beta}+
f_{\beta,\xi}\zeta_{\beta,\xi}.
\end{equation}
Then $\omega(\xi,\eta,\beta)$ 
is a continuous bounded function on $G\times G$.
Since $f_{\xi,\eta}=-f_{\eta,\xi}$ and 
$\zeta_{\xi,\eta}=\zeta_{\eta,\xi}$ for all
$\xi\not=\eta\in \Lambda$, 
we have $\omega\circ \sigma=({\rm
sgn}(\sigma))\omega$ for every permutation $\sigma$ of the three
variables. Since the functions $f_{\xi,\eta},\zeta_{\xi,\eta}$
are restrictions to the leaves of the foliation ${\cal F}$
of globally continuous bounded functions on $\tilde V$, 
the assignment $(\xi,\eta,\beta)\in T\to 
\omega(\xi,\eta,\beta)\in C^0(G\times G)$ is continuous with respect to the
compact open topology on $C^0(G\times G)$. Moreover, it, is
equivariant
with respect to the natural action of $G$ on $T$ and on $C^0(G\times G)$.
This means that $\omega$ is a continuous bounded cocycle 
for the action of $G$ on $\Lambda$ with
values in $C^0(G\times G)$.

{\sl Step 3:}

Recall the choices of the
constants $0<\kappa_0<\kappa<\kappa_1$ from Step 1 above. 
We may assume that $\kappa >1$. 
In this technical step we obtain some control on the 
functions $f_{\xi,\eta}\zeta_{\xi,\eta}$ for a pair
of distinct points $\xi\not=\eta \in \Lambda$.

By hyperbolicity of $X$ and the choice of $\kappa_0$, for every
triple $(\xi,\eta,\beta)$ of pairwise distinct points in $\partial X$ there
is a point $y_0\in X$ which is contained in the
$\kappa_0$-neighborhood of every side of a geodesic triangle with
vertices $\xi,\eta,\beta$. Let $\gamma:\mathbb{R}\to X$ be a geodesic
connecting $\xi$ to $\eta$ with $d(y_0,\gamma(0))\leq \kappa_0.$ Then
the Gromov product $(\xi,\beta)_{\gamma(t)}$ based at $\gamma(t)$
between $\xi,\beta$ satisfies 
\begin{equation}
(\xi,\beta)_{\gamma(t)}\geq t-c_0\text{ for all }
t\geq 0\notag\end{equation}
with a constant $c_0>0$ only depending on $X$
(see \cite{BH}). Thus
by the properties (\ref{dist}) and (\ref{distequi}) of the distance
functions $\delta_x$, 
there is a number $r_0>0$ (depending on $\kappa_1$ and
the hyperbolicity constant $\kappa_0$ for $X$) 
such that if $t\geq 0$ and if $y\in X$
satisfies $d(\gamma(t),y)<3\kappa_1$ then 
$\delta_y(\xi,\beta)\leq
r_0e^{-\chi t}$
where $\chi>0$ is as in inequality (\ref{dist}).
The triangle inequality for $\delta_y$ then yields 
\begin{equation}\label{deltadifference}
\vert\delta_y(\xi,\eta)-\delta_y(\eta,\beta)\vert \leq r_0e^{-\chi t}.
\end{equation} 

This implies the following. First, recall that
the auxiliary function $\phi_0$ which we used in the definition 
(\ref{phi}) of the
function $\phi$ is smooth and constant outside a compact set and 
hence it is uniformly
Lipschitz continuous. Therefore by the
estimate (\ref{deltadifference}) there is a constant
$r_1>r_0$ such that
\begin{equation}\label{phi2}
\vert\phi_0(\delta_y(\xi,\eta))-\phi_0 (\delta_y(\eta,\beta))\vert \leq
r_1e^{-\chi t}
\end{equation} whenever $d(y,\gamma(t)) \leq 3\kappa_1$
for some $t\geq 0$.

Now let  
$0\leq t$ and let $u,h\in G$ be such that
\[d(ux_0,\gamma(t))<3\kappa_1,d(hx_0,\gamma(t))<3\kappa_1.\] 
Since the function 
$\phi_0$ assumes values in $[0,1]$, we obtain from the
definition (\ref{phi}) of the function $\phi$ and the estimate 
(\ref{phi2}) that
\begin{align}\label{phi3}
\vert \phi(\xi,\eta,ux_0,hx_0)& -\phi(\beta,\eta,ux_0,hx_0)\vert \\ \leq
\vert (\phi_0(\delta_{ux_0}(\xi,\eta))& -
\phi_0(\delta_{ux_0}(\beta,\eta)))\phi_0(\delta_{hx_0}(\xi,\eta))\vert \notag\\
+\vert \phi_0(\delta_{ux_0}(\beta,\eta))(\phi_0(\delta_{hx_0}(\xi,\eta)) &
-\phi_0(\delta_{hx_0}(\beta,\eta)))\vert \leq 2r_1e^{-\chi t}.\notag
\end{align} 
Similarly, 
the function $\psi$ used in the definition (\ref{zeta}) of the
function $\zeta$ assumes values in $[0,1]$ and hence
we conclude from (\ref{phi3}) that also
\begin{equation}\label{zeta2}
\vert \zeta(\xi,\eta,ux_0,hx_0)-\zeta(\beta,\eta,ux_0,hx_0)\vert \leq
2r_1e^{-\chi t}.
\end{equation}
Moreover, by the definition (\ref{productmetric2}) 
of the distance function $\hat d$ on $V$ and by
the inequality (\ref{tilded}) for the distance
function $\tilde d$ on $\partial X\times X$, 
we have 
\begin{equation}\label{hatdestimate}
\hat d\bigl((\xi,\eta,hx_0,ux_0),(\beta,\eta,hx_0,ux_0)\bigr)
\leq\frac{1}{2}(\delta_{hx_0}(\xi,\beta)+\delta_{ux_0}(\xi,\beta))
\leq r_0e^{-\chi t}.
\end{equation}

The function $\hat f:V\to \mathbb{R}$ 
constructed in Step 2 from a function $f\in {\cal H}$ 
is H\"older continuous and $\iota$-anti-invariant. Therefore
by the estimate (\ref{hatdestimate})
there are numbers $\alpha>0,r_2>r_1$ only depending on the H\"older
norm for $f$ with the following property.
Let $0\leq t$ and let $u,h\in G$ be such that
$d(ux_0,\gamma(t))<3\kappa_1,d(hx_0,\gamma(t))<3\kappa_1$; 
then\begin{equation}\label{hatf}
\vert \hat f(\xi,\eta,ux_0,hx_0)+
\hat f(\eta,\beta,ux_0,hx_0)\vert \leq r_2e^{-\chi \alpha t}.
\end{equation}
The functions $f$ and $\zeta$ are bounded in absolute value
by a universal constant. Hence using a calculation as
in (\ref{phi3}) above, from 
the definition of the functions $f_{\xi,\eta}$ and
$f_{\eta,\beta}$ and from the estimates (\ref{zeta2}) and (\ref{hatf})
we obtain the existence of a constant
$r>r_2$ (depending on the H\"older norm of $f$) 
such that
\begin{equation}\label{asymptotic}
\vert (f_{\xi,\eta}\zeta_{\xi,\eta}+
f_{\eta,\beta}\zeta_{\eta,\beta})(u,h)\vert \leq re^{-\chi \alpha s}.
\end{equation}

{\sl Step 4:}

Let $\nu=\mu\times \mu$ be the left invariant product measure on
$G\times G$. Our goal is to show that for every $p\in (0,\infty)$, 
the cocycle
$\omega$ defined in equation (\ref{omega}) above 
is in fact a bounded cocycle with values in
$L^p(G\times G,\nu)$. 
For this we show that for every $(\xi,\eta,\beta)\in T$ 
the function $\omega(\xi,\eta,\beta)$ on
$G\times G$ is contained in a fixed bounded subset of 
$L^p(G\times G,\nu)$ and that
moreover the assignment 
$(\xi,\eta,\beta)\to \omega(\xi,\eta,\beta)\in L^p(G\times G,\nu)$ is
continuous.

For a
subset $C$ of $X$ write 
\[C_{G,2\kappa_1}=\{(u,h)\in G\times G\mid
ux_0\in C,d(ux_0,hx_0)\leq 2\kappa_1\}\] 
and for $r>0$ let $N(C,r)$
be the $r$-neighborhood of $C$ in $X$. 
We claim that 
there is a number $m>0$ such that for every subset
$C$ of $X$ of diameter at most one the 
$\nu$-mass of the set $N(C,2\kappa)_{G,2\kappa_1}$ 
is at most $m$. 
Namely, the subset 
\[D=\{(u,h)\in G\times G\mid d(ux_0,x_0)\leq 6\kappa,
d(ux_0,hx_0)\leq 2\kappa_1\}\] 
of $G\times G$ is compact and hence
its $\nu$-mass is finite, say this mass equals $m>0$. 
On the other hand, if $C\subset X$ is a set of
diameter at most one and if there is some 
$g\in G$ such that $gx_0\in N(C,2\kappa)$ then
any pair $(u,h)\in N(C,2\kappa)_{G,2\kappa_1}$ is contained
in $gD$. Our claim now
follows from the fact that $\nu$ is invariant under left
translation.

As in Step 3 above, let $(\xi,\eta,\beta)$ a triple of
pairwise distinct points in $\Lambda$ and let
$y_0\in X$ be a point which is contained in the
$\kappa_0$-neighborhood of every side of a geodesic
triangle in $X$ with vertices $\xi,\eta,\zeta$.
Let $\gamma:\mathbb{R}\to X$ be a geodesic connecting
$\xi$ to $\eta$ with $d(y_0,\gamma(0))\leq \kappa_0$.
Also, let $\rho:\mathbb{R}\to X$ 
be a geodesic connecting $\beta$ to $\eta$ which is
parametrized in such a way that $d(y_0,\rho(0))\leq \kappa_0$. 
Then $\gamma[0,\infty),\rho[0,\infty)$ are two sides of 
a geodesic triangle in $X$ with vertices $\gamma(0),\rho(0),\eta$.
Since $\kappa_0$ is a hyperbolicity constant for $X$
and $d(\gamma(0),\rho(0))\leq 2\kappa_0$, the
ray $\gamma[3\kappa_0,\infty)$ is contained in the
$\kappa_0$-tubular neighborhood of $\rho[0,\infty)$, and the
ray $\rho[3\kappa_0,\infty)$ is contained in the 
$\kappa_0$-tubular neighborhood of $\gamma[0,\infty)$. Thus
the Hausdorff distance between
the geodesic rays $\gamma[0,\infty)$ and $\rho[0,\infty)$ is at
most $6\kappa_0\leq \kappa$. 
In particular, the $\kappa$-neighborhood of $\rho[0,\infty)$
is contained in the $2\kappa$-neighborhood of
$\gamma[0,\infty)$.

Recall that
$\phi_0(\delta_{ux_0}(\xi,\eta))\not=0$ only if $ux_0$ is contained in
the $\kappa$-neighborhood of $\gamma$.
Moreover, we have $\zeta(\xi,\eta,ux_0,hx_0)\not=0$ only
if $d(ux_0,hx_0)\leq 2\kappa_1$. This implies that the
support of the function $f_{\xi,\eta}\zeta_{\xi,\eta}$ is contained
in $N(\gamma(\mathbb{R}),\kappa)_{G,2\kappa_1}$ and similarly
for the functions $f_{\eta,\beta}\zeta_{\eta,\beta},
f_{\beta,\xi}\zeta_{\beta,\xi}$.

As a consequence, there is an open subset $U$ of $G\times G$
with compact closure 
and there is a number $T>0$ only depending on $\kappa_1$
and the hyperbolicity constant for $X$ with the following
property. The support of the function $\omega$ defined in
(\ref{omega}) above is the disjoint union of the 
three sets \[N(\gamma[T,\infty),2\kappa)_{G,2\kappa_1},
N(\gamma(-\infty,-T],2\kappa)_{G,2\kappa_1},
N(\rho(-\infty,-T],2\kappa)_{G,2\kappa_1}\] 
with $U$.
Moreover, the restriction of $\omega$ to
$N(\gamma[T,\infty),2\kappa)_{G,2\kappa_1}$ coincides with
the restriction of the function
$f_{\xi,\eta}\zeta_{\xi,\eta}+f_{\eta,\beta}\zeta_{\eta,\beta}$
and similarly for the other two sets different from $U$ in the
above decomposition of the support of $\omega$. Thus to show
that $\omega$ is contained in $L^p(G\times G)$ it is enough to
show that there is constant $c_p>0$ only depending on $p$ and 
the H\"older norm of $f$ such that
\begin{equation}\label{integralestimate1}
\int_{N(\gamma[T,\infty),\kappa_1)_{G,\kappa_2}}\vert
f_{\xi,\eta}\zeta_{\xi,\eta}+ 
f_{\eta,\beta}\zeta_{\eta,\beta}\vert^p d\nu < c_p. \notag
\end{equation}

However, this is immediate from the 
estimate (\ref{asymptotic}) 
together with the control on
the $\nu$-mass of subsets of   
$N(\gamma[T,\infty),2\kappa)_{G,2\kappa_1}$ $(k>0)$.
Namely, we showed that for every integer $k\geq 0$ the
$\nu$-mass of the set
$N(\gamma[T+k,T+k+1],2\kappa)_{G,2\kappa_1}$ is bounded from
above by a universal constant $m>0$. 
Moreover, for every $p\geq 1$ the value of the function
$\vert f_{\xi,\eta}\zeta_{\xi,\eta}+f_{\beta,\eta}\zeta_{\beta,\eta}\vert^p $
on this set does not exceed $r^pe^{-p\chi\alpha(T+k)}$.
Thus the inequality \label{integralestimate} holds true with
$c_p=mr^p\sum_{k=0}^\infty e^{-p\chi\alpha(T+k)}$.

Since the function $\tilde f$ on $\tilde V$ is globally continuous,
the same consideration also
shows that $\omega(\xi,\eta,\beta)\in L^p(G\times G,\nu)$
depends continuously on $(\xi,\eta,\beta)$. 
Namely, let $((\xi_i,\zeta_i,\eta_i))\subset T$ be a sequence
of triples of pairwise distinct points converging to a 
triple $(\xi,\eta,\beta)\in T$. By the above consideration,
for every $\epsilon >0$ there is a compact subset $A$ of $G\times G$
such that $\int_{G\times G-A}\vert \omega(\xi_i,\eta_i,\beta_i)\vert^p
d\nu\leq \epsilon$ for all sufficiently large $i>0$ and that
the same holds true for $\omega(\xi,\eta,\zeta)$. 
Let $\chi_A$ be the characteristic function of $A$.
By continuity of the
function $\tilde f\tilde \zeta$ on $\tilde V$ and compactness,
the functions $\chi_A\omega(\xi_i,\eta_i,\beta_i)$ converge
as $i\to \infty$ in $L^p(G\times G,\nu)$ to $\chi_A\omega(\xi,\eta,\zeta)$.
Since $\epsilon >0$ was arbitrary, the required continuity follows.

Moreover, the assignment
$(\xi,\eta,\beta)\to \omega(\xi,\eta,\beta)$ 
is equivariant under the action of $G$
on the space $T$ of triples of pairwise distinct points in
$\Lambda$ and on $L^p(G\times G,\nu)$ and satisfies the
cocycle equality (\ref{coc}).
In other words, $\omega$ defines a continuous 
$L^p(G\times G,\nu)$-valued bounded cocycle for the action of $G$
on $\Lambda$ as required.

{\sl Step 5:}

Let $g\in G$ be a hyperbolic isometry, let $a\not=b\in \partial X$
be the attracting and repelling 
fixed point for the action of $g$ on $\partial X$, respectively,
and let $\xi\in \Lambda-\{a,b\}$. We have to show
that we can find a cocycle $\omega$ as in (\ref{omega}) above with
$\omega(a,b,\xi)\not=0$.

For this let $\gamma$ be
a geodesic in $X$ connecting $b$ to $a$ and choose the basepoint
$x_0$ for the above construction 
on $\gamma$. Let $\kappa>\kappa_0$ be as above. Then
the orbit of $x_0$ under the infinite cyclic subgroup of 
$G$ generated by $g$ is contained in the $\kappa_0$-neighborhood
of $\gamma$ (compare the discussion in the proof of 
Lemma \ref{fixedpoint}).
Since $g^jx_0\to a,g^{-j}x_0\to b$ $(j\to \infty)$, there
are numbers 
$k<\ell$ such that the
$\kappa$-neighborhood of a geodesic connecting $a$ to $\xi$ and of
a geodesic connecting $b$ to $\xi$ contains at most one of the
points $g^kx_0,g^\ell x_0$ and that moreover the distance
between $g^kx_0,g^\ell x_0$ is at least $4\kappa_0$. Choose 
$\kappa_1>2d(g^kx_0,g^\ell x_0)$.

Using this constant $\kappa_1$ for the
definition of the function $\zeta$ in equation (\ref{zeta}),
by the definition of $\phi_0$ we have
$\zeta(a,b,g^kx_0,g^\ell x_0)>0$ and 
$\zeta(a,\xi,g^kx_0,g^\ell x_0)=\zeta(b,\xi,g^kx_0,g^\ell x_0)=0$.
Thus the pair of points
$(g^k,g^\ell)\in G\times G$ is contained in the support of the function
$\zeta_{a,b}$ but not in the support of any of the two functions
$\zeta_{a,\xi},\zeta_{b,\xi}$. Moreover, 
by Lemma \ref{fixedpoint}, there is no $h\in G$ with
$hg^kx_0=g^kx_0,hg^\ell x_0=g^\ell x_0$ and $h(a)=b,h(b)=a$.
Therefore the $G$-orbit of
$(a,b,g^kx_0,g^\ell x_0)\in V$ 
does not contain the point $(b,a,g^kx_0,g^\ell x_0)$.
This means that the 
projection of $(a,b,g^kx_0,g^\ell x_0)$ 
into $W$ is not fixed by the involution $\iota$.

As a consequence, we can find a function $f\in {\cal
H}$ (for a suitable choice of a support ball $B$)
whose lift $\tilde f$ to $\tilde V$ 
does not vanish at
$(a,b,g^k,g^\ell)$. By the choice of $\kappa_1$,
this means that $f_{a,b}\zeta_{a,b}(g^k,g^\ell)\not=0$
and $f_{b,\xi}\zeta_{b,\xi}(g^k,g^\ell)=
f_{\xi,a}\zeta_{\xi,a}(g^k,g^\ell)=0$. In other words, the
cocycle $\omega$ constructed as above
from $f$ does not vanish at $(a,b,\xi)$. This shows the theorem.
\end{proof}

{\bf Remark:} The construction in the proof of Theorem \ref{cocycle}
yields in fact for every triple $(a,b,\xi)\in \Lambda^3$ such that
$(a,b)$ is the pair of fixed points of a hyperbolic element of $G$
an infinite dimensional space of continuous 
$L^p(G\times G,\mu\times \mu)$-valued cocycles 
which do not vanish at $(a,b,\xi)$.

\section{Bounded cohomology}

In this section we use Theorem \ref{cocycle} to construct 
nontrivial second bounded cohomology classes for 
closed non-elementary sugbroups of the isometry
group ${\rm Iso}(X)$
of a proper hyperbolic geodesic metric space $X$.

Every locally compact $\sigma$-compact topological group $G$ admits a
\emph{strong boundary} \cite{K03} which is a standard Borel
$G$-space $B$ with a quasi-invariant ergodic probability measure
$\lambda$ such that the action of $G$ on $(B,\lambda)$ is amenable
and doubly ergodic with respect to any separable Banach module
for $G$ (\cite{M}, Definition 11.1.1). The following lemma is
basically contained in \cite{MS04} and follows from the arguments
of Zimmer (see \cite{Z}). For later use we formulate it more
generally for \emph{cocycles}. 

Namely, let
$G$ be a locally compact $\sigma$-compact group which admits a
measure preserving ergodic action on a standard probability space
$(S,\nu)$. An ${\rm Iso}(X)$-valued \emph{cocycle} for this action
is a measurable map $\alpha:G\times S\to {\rm Iso}(X)$ such that
\begin{equation}
\alpha(gh,x)=\alpha(g,hx)\alpha(h,x) \notag
\end{equation}
for all $g,h\in G$ and
almost all $x\in S$. The cocycle $\alpha$ is \emph{cohomologous}
to a cocycle $\beta$ if there is a measurable map $\psi:S\to {\rm
Iso}(X)$ such that $\psi(gx)\alpha(g,x)=\beta(g,x)\psi(x)$ for
almost every $x\in S$, all $g\in G$. We have.

\begin{lemma}\label{reduction}
Let $X$ be a proper hyperbolic geodesic
metric space with isometry group ${\rm Iso}(X)$. Let $G$ be a
locally compact $\sigma$-compact group which admits a measure
preserving ergodic action on a standard probability space
$(S,\nu)$ and let $\alpha:G\times S\to {\rm Iso}(X)$ be a cocycle.
Assume that $\alpha(G\times S)$ is contained in a closed
non-elementary subgroup $H$ of ${\rm Iso}(X)$ and that $\alpha$ is
not cohomologous to a cocycle into 
an elementary subgroup of $H$. Let
$\Lambda$ be the limit set of $H$ and let $(B,\lambda)$ be a
strong boundary for $G$; then there is an $\alpha$-equivariant
measurable map $\Psi:B\times S\to \Lambda$.
\end{lemma}
\begin{proof}
Let $X$ be a proper hyperbolic geodesic metric space
and let $G$ be a locally compact $\sigma$-compact group with a
measure preserving ergodic action on a standard probability space
$(S,\nu)$. Let $\alpha:G\times S\to
{\rm Iso}(X)$ be a cocycle. Assume that  
the closure $H$ of $\alpha(G\times S)$ in
${\rm Iso}(X)$ is non-elementary, with limit set $\Lambda$, and
that $\alpha$ is not cohomologous to a cocycle into an elementary
subgroup of $H$. Let $(B,\lambda)$ be a strong boundary
for $G$. Since the action of $G$ on $B$ is amenable, there is an
$\alpha$-equivariant Furstenberg map $f:B\times S\to {\cal
P}(\Lambda)$ where ${\cal P}(\Lambda)$ denotes the space of Borel
probability measures on $\Lambda$ \cite{Z}.

Let ${\cal P}_{\geq 3}(\Lambda)$ be the space of probability
measures on $\Lambda$ whose support contains at least 3 points. By
Corollary 5.3 of \cite{A96}, the $H$-action on ${\cal P}_{\geq
3}(\Lambda)$ is tame with compact point stabilizers. Moreover,
since $B$ is a strong boundary for $G$, the action of $G$ on
$B\times S$ is ergodic. Thus if
the image under $f$ of a set of positive measure in $B\times S$ is
contained in ${\cal P}_{\geq 3}(\Lambda)$ then by ergodicity we
have $f(B\times S)\subset {\cal P}_{\geq 3}(\Lambda)$ (almost
everywhere). Moreover by the results of Zimmer \cite{Z},
the cocycle $\alpha$ is cohomologous to a cocycle into a
compact subgroup of $H$. This is a contradiction to
our assumption on $\alpha$.

As a consequence, we have
$f(B\times S)\subset {\cal P}_{\leq 2}(\Lambda)$
where ${\cal P}_{\leq 2}(\Lambda)$ is the set of measures whose
support contains at most 2 points. To show that the
measures in $f(B\times S)$
are in fact supported in a single point of $\Lambda$ we
argue as in the proof of Lemma 3.4 of \cite{MS04}. Namely, the
group $H$ is locally compact and $\sigma$-compact and its action
on the space of triples of pairwise distinct points in $\Lambda$
is proper (see \cite{A96}). Thus the assumptions in
Lemma 23 of \cite{MMS04} are satisfied. We can then use Lemma 23
of \cite{MMS04} as in the proof of Lemma 3.4 
of \cite{MS04} to conclude that the image of $f$ is in fact
contained in the set of Dirac masses on $\Lambda$, i.e. there is
an $\alpha$-equivariant map $B\times S\to \Lambda$ as claimed. 
\end{proof}

We use Lemma \ref{reduction} and Theorem \ref{cocycle} to deduce.

\begin{corollary} \label{nonzero}
Let $G$ be a locally compact
$\sigma$-compact group and let $\rho$ be a homomorphism of $G$
into the isometry group ${\rm Iso}(X)$ 
of a proper hyperbolic geodesic metric
space $X$. Let $H$ be the closure of 
$\rho(G)$ in ${\rm Iso}(X)$ and let $\mu$ be Haar measure on $H$.
If $H$ is
non-elementary, then $H_{cb}^2(G,L^p(H,\mu))\not=\{0\}$
for every $p\in (1,\infty)$.
\end{corollary}
\begin{proof} 
Let $\rho:G\to {\rm Iso}(X)$ be a homomorphism and let $H$ be the
closure of $\rho(G)$ in ${\rm Iso}(X)$. Assume that $H$ is
non-elementary and let $\Lambda\subset \partial X$ be the limit
set of $H$. By Lemma \ref{reduction},
applied to the homomorphism $\rho$ viewed as a cocycle
for the trivial action of $G$ on a point, 
there is a measurable 
$\rho$-equivariant map
$\phi$ from a strong boundary $(B,\lambda)$
of $G$ into $\Lambda$. Since $\rho(G)$ is dense in $H$, the
set $\Lambda$ is also the limit set of $\rho(G)$ and therefore
the $\rho(G)$-orbit of every point in $\Lambda$ is
dense in $\Lambda$ (8.27 in \cite{GH}).
By equivariance, the image under $\phi$ of $\lambda$-almost 
every $G$-orbit on $B$ is
dense in $\Lambda$. In particular, the support of the
measure class $\phi_*\lambda$
is all of $\Lambda$.

Let $\mu$ be the Haar measure of $H$.
By Theorem \ref{cocycle}, for every $p\in (1,\infty)$ 
there is a nontrivial bounded continuous
$L^p(H\times H,\mu\times \mu)$-valued cocycle $\omega$ on the
space of triples of pairwise distinct points in $\Lambda$.
Then the
$L^p(H\times H,\mu\times \mu)$-valued
$\lambda$-measurable bounded cocycle $\omega\circ \phi^3$
on $B\times B\times B$ is non-trivial on a set of positive
measure. Since $B$ is a strong boundary for $G$, this cocycle then
defines a \emph{non-trivial} class in $H_{cb}^2(G,L^p(H\times
H,\mu\times \mu))$ (see \cite{M}). On the other hand, the isometric
$G$-representation space $L^p(H\times H,\mu\times \mu)$ is a
direct integral of copies of the isometric $G$-representation space
$L^p(H,\mu)$ and therefore by Corollary 2.7 of \cite{MS04} and
Corollary 3.4 of \cite{MS05}, if $H_{cb}^2(G,L^p(H,\mu))=\{0\}$
then also $H_{cb}^2(G,L^p(H\times H,\mu\times \mu))=\{0\}$. This
implies the corollary. \end{proof}

Corollary \ref{nonzero} applied to a closed non-elementary
subgroup $G<{\rm Iso}(X)$ shows that $H_{cb}^2(G,L^p(G,\mu))\not=\{0\}$
for every $p\in (1,\infty)$.

\section{Structure of the isometry group}

In this section we analyze the structure of the isometry group of
a proper hyperbolic geodesic metric space $X$ using the results of
Section 2 and Section 3. 
For this recall that an \emph{extension} of a locally
compact group $H$ by a topological group $K$ is a locally compact
group $G$ which contains $K$ as a normal subgroup and such that
$H=G/K$ as topological groups. The following easy observation of
Monod and Shalom \cite{MS04} describes non-elementary closed
subgroups of ${\rm Iso}(X)$ which are extensions by amenable
groups.

\begin{lemma}\label{amenable} 
Let $X$ be a proper hyperbolic geodesic
metric space and let $G$ be a closed non-elementary
subgroup of ${\rm Iso}(X)$.
Then a maximal normal amenable subgroup
of $G$ is compact.
\end{lemma}
\begin{proof}
Let $G$ be a closed non-elementary
subgroup of ${\rm Iso}(X)$ and
let $H\lhd G$ be a maximal normal amenable subgroup. Since $H$ is
amenable and the Gromov boundary $\partial X$ of $X$ is compact,
there is an $H$-invariant probability measure $\mu$ on $\partial
X$. Thus $H$ is compact if the support $A$ of $\mu$ contains at
least three points (see the proof of Theorem 21 in \cite{MMS04}).

Now assume that $A$ contains at most two points.
Since $G$ is non-elementary by assumption, 
the limit set of $G$ is the smallest $G$-invariant 
closed subset of $\partial X$ (8.27 of \cite{GH}). 
Thus the set $A$ is not invariant 
under the action of $G$. Let $g\in G$ be
such that $A\cup gA\cup g^2A$ contains at least 3 points. Since
$H$ is a normal subgroup of $G$ and
$A$ is $H$-invariant, the set $A\cup gA\cup g^2A$ is
preserved by $H$ and hence as before, the group $H$ is necessarily
compact. \end{proof}

Let $G$ be a locally compact $\sigma$-compact topological group
and let $N\lhd G$ be a compact normal subgroup. For an isometric
representation $\rho$ of $G$ into a separable Banach space $E$
let $E^N$ be the closed subspace of $E$ of all $N$-invariant
vectors. Then $\rho$ induces a representation of $G/N$ into $E^N$.
Since compact groups are amenable the following observation is a
special case of Corollary 8.5.2 of \cite{M}.

\begin{lemma}\label{representation} 
Let $G$ be an extension of a topological group $H$ by a
compact group $N$ and let $\rho$ be an isometric representation of
$G$ into a separable Banach space $E$. Then the projection
$P:G\to H$ induces an isomorphism of $H_{cb}(H,E^N)$ onto
$H_{cb}(G,E)$.
\end{lemma}

We use Corollary \ref{nonzero}, 
Lemma \ref{amenable} and Lemma \ref{representation} to complete the
proof of Theorem 1 from the introduction.

\begin{proposition}\label{firstpart} 
Let $X$ be a proper hyperbolic
geodesic metric space and let $G< {\rm Iso}(X)$ be a closed
subgroup. Then one of the following three
possibilities holds.
\begin{enumerate}
\item $G$ is elementary.
\item  $G$ contains an open subgroup $G^\prime$ of
finite index which
is a compact extension of a simple Lie group of rank 1.
\item $G$ is a compact extension of a totally disconnected
group.
\end{enumerate}
\end{proposition}
\begin{proof}
Let $G$ be a closed subgroup of the isometry group
${\rm Iso}(X)$ of a proper hyperbolic geodesic metric space $X$.
Then $G$ is locally compact. Assume that $G$ is non-elementary;
then by Lemma \ref{amenable}, the maximal normal amenable
subgroup $H$ of $G$ is compact, and the quotient $V=G/H$ is a
locally compact $\sigma$-compact topological group. Denote by
$\pi:G\to G/H=V$ the canonical projection.

By the solution to Hilbert's fifth problem (see Theorem 11.3.4 in
\cite{M}), after possibly replacing $G$ by an open subgroup of
finite index (which we denote again by $G$ for simplicity), the
quotient $V=G/H$ splits as a direct product $V=V_0\times Q$ where
$V_0$ is a semi-simple connected Lie group with finite center and
without compact factors and $Q$ is totally disconnected. If $V_0$
is trivial then $G$ is a compact extension of a totally
disconnected group.

Now assume that $V_0$ is nontrivial. Then $V_0$ is not compact and
the limit set $\Lambda$ of $\pi^{-1}(V_0)< G$ is nontrivial. Since
$Q$ commutes with $V_0$, the group $\pi^{-1}(Q)<G$ commutes with
$\pi^{-1}(V_0)<G$ up to the compact normal subgroup $H$ and hence
$\pi^{-1}(Q)$ acts trivially on $\Lambda$. Namely, let $x\in X$
and let $(g_i)\subset \pi^{-1}(V_0)$ be a sequence such that
$(g_ix)$ converges to $\xi\in \Lambda$. If $h\in Q$ is arbitrary
then for every $i>0$ there is some $h_i^\prime\in H$ with
$hg_i=g_ihh_i^\prime$ for all $i$.
Now $H$ is compact and therefore the sequence $(g_ihh_i^\prime x)$
converges as $i\to \infty$ to $\xi$. This implies that
$(hg_i x)$ converges
to $h\xi=\xi$ $(i\to\infty)$. 
Thus if $\Lambda$ consists of one or
two points then the action of $G$ on $\partial X$ fixes one or two
points. Since $G$ is non-elementary by assumption, this is
impossible. 

As a consequence, $\Lambda$
contains at least three points. Then $\pi^{-1}(Q)$ fixes at least
three points in $\partial X$ and hence $\pi^{-1}(Q)< G$ is
necessarily a \emph{compact} normal subgroup of $G$ containing
$H$. Since by assumption $H$ is a maximal compact normal subgroup
of $G$ we conclude that $\pi^{-1}(Q)=H$ and that $Q$ is trivial.
In other words, if $G$ is non-elementary then up to passing to an
open subgroup of finite index, either $G$ is a compact extension
of a totally disconnected group
or $G$ is a compact extension of a
semi-simple Lie group $V_0$ with finite center and without compact
factors.

We are left with showing that if $G$ is a compact extension of a
semi-simple Lie group $V_0$ with finite center and without compact
factors then $V_0$ is simple and of rank 1.
For this note from Corollary \ref{nonzero} and 
Lemma \ref{representation} that
$H_{cb}^2(V_0,L^2(V_0,\mu))\not=\{0\}$. On the other hand, if $V_0$
is a \emph{simple} Lie group of non-compact type and rank at least
2 then Theorem 1.4 of \cite{MS04} shows that
$H_{cb}^2(V_0,L^2(V_0,\mu))=\{0\}$. The case that 
$V_0=V_1\times V_2$ is
a non-trivial product of semi-simple Lie groups $V_1,V_2$
of non-compact type can be ruled out in the same way. Namely, 
in this case the representation space $L^2(V_0,\mu)$
does not admit any nontrivial $V_i$-invariant vectors 
$(i=1,2)$ and therefore
$H_{cb}^2(V_0,L^2(V_0,\mu))=\{0\}$ by the Burger-Monod super-rigidity
result for cohomology \cite{BM02} (see Theorem 14.2.2 in
\cite{M}). Together we conclude that necessarily $V_0$ is a simple
Lie group of rank one. This finishes the proof of the proposition.
\end{proof}

\section{Super-rigidity of cocycles}

Let again $X$ be a proper hyperbolic geodesic metric space. 
The goal of this section is to study ${\rm Iso}(X)$-valued
cocycles and derive Theorem 3 from the introduction.

Let $S$ be a standard Borel space and let
$\nu$ be a Borel probability
measure on $S$. Let $\Gamma$ be any countable group 
which admits a measure preserving ergodic action on
$(S,\nu)$. This action then defines a natural
continuous unitary representation of $\Gamma$ into
the Hilbert space $L^2(S,\nu)$ of square integrable
functions on $S$.
Let $\alpha:\Gamma\times S\to {\rm Iso}(X)$ be a cocycle,
i.e. $\alpha:G\times S\to {\rm Iso}(X)$ is a Borel
map which satisfies $\alpha(gh,x)=\alpha(g,hx)\alpha(h,x)$
for all $g,h\in \Gamma$ and $\nu$-almost every $x\in S$.

If $\Gamma$ is a lattice in a locally compact $\sigma$-compact
topological group $G$ then 
$\Gamma$ admits a measure preserving action on the product 
space $G\times S$. Since $\Gamma<G$ is a lattice, the
quotient space $(G\times S)/\Gamma$ can be
viewed as a bundle over $G/\Gamma$ with fibre $S$.
If $\Omega\subset G$ is a 
finite measure Borel fundamental domain
for the action of $\Gamma$ on $G$ then $\Omega\times S\subset
G\times S$ is a finite measure Borel fundamental domain for
the action of $\Gamma$ on $G\times S$. Thus if
we denote by $\mu$ a Haar measure on $G$ then 
up to normalization, the
product measure $\mu\times \nu$ projects to a probability
measure $\lambda$ on $(G\times S)/\Gamma$. The action of $G$
on $G\times S$ by left translation commutes with the
action of $\Gamma$ and hence it projects to an action on
$(G\times S)/\Gamma$ preserving the measure $\lambda$ 
(see p.75 of \cite{Z}).   

Recall that the action of $\Gamma$ on $(S,\nu)$ is mildly mixing
if $\lim\inf_{g\to \infty}\nu(A\Delta gA)\not=0$ for every
Borel set $A\subset S$ with $0<\nu(A)<1$.
The following lemma is due to Schmidt and Walters \cite{SW82}.

\begin{lemma}\label{mildmixing} 
Let $\Gamma$ be an irreducible lattice in a
product $G=G_1\times G_2$ of two semi-simple non-compact Lie
groups $G_1,G_2$ with finite center.
If the action of $\Gamma$ on $(S,\nu)$ is mildly mixing
then the induced action of $G_1$ on $((G\times S)/\Gamma,\lambda)$
is ergodic.
\end{lemma}

The following proposition completes the proof of
Theorem 3 from the introduction and follows as in
\cite{MS04} from Lemma \ref{mildmixing}, the
rigidity results for bounded cohomology
of Burger and Monod \cite{BM99,BM02} and
the work of Zimmer \cite{Z}.

\begin{proposition} \label{cocyclerig}
Let $G$ be a semi-simple Lie group with
finite center, no compact factors and of rank at least $2$. Let
$\Gamma<G$ be an irreducible lattice and let $(S,\nu)$ be a 
mildly mixing
$\Gamma$-space. Let $X$ be a proper hyperbolic geodesic metric
space and let $\alpha:\Gamma\times S\to {\rm Iso}(X)$ be a
cocycle. Then either $\alpha$ is cohomologous to a cocycle into an
elementary subgroup of ${\rm Iso}(X)$ or there is a closed
subgroup $H$ of ${\rm Iso}(X)$ which is a compact extension of a
simple Lie group $L$ of rank one and there is a surjective
homomorphism $G\to L$.
\end{proposition}
\begin{proof}
Let $G$ be a semi-simple Lie group of non-compact
type with finite center and of rank at least 2 and let $\Gamma$ be
an irreducible lattice in $G$. Let $(S,\nu)$ be a mildly mixing
$\Gamma$-space with invariant Borel probability measure $\nu$. Let
$\alpha:\Gamma\times S\to {\rm Iso}(X)$ be a cocycle into the
isometry group of a proper hyperbolic geodesic metric space $X$.
Assume that $\alpha$ is not cohomologous to a cocycle
into an elementary subgroup of ${\rm Iso}(X)$.
Let $H<{\rm Iso}(X)$ be a closed subgroup with limit
set $\Lambda$ which contains the image of $\alpha$.

Let $\Omega\subset G$ be a Borel fundamental domain for
the action of $\Gamma$ on $G$. Let $\lambda$ be the $G$-invariant
Borel probability measure on $(G\times S)/\Gamma = \Omega\times S$. 
We obtain
a $\lambda$-measurable function $\beta:G\times (G\times S)/\Gamma\to 
H$ as follows. For $z\in \Omega$ and $g\in G$
let $\eta(g,z)\in \Gamma$ be the unique
element such that $gz\in \eta(g,z)\Omega$ and define
$\beta(g,(z,\sigma))=\alpha(\eta(g,z),\sigma)$. By construction,
$\beta$ satisfies the cocycle equation for the action of $G$ on
$(G\times S)/\Gamma$. 

Let $(B,\rho)$ be a strong boundary for $G$; we may assume that $B$ is
also a strong boundary for $\Gamma$. By Lemma \ref{mildmixing}, 
the action
of $G$ on $(G\times S)/\Gamma$ is ergodic and hence 
by Lemma \ref{reduction}, there is 
a measurable $\beta$-equivariant
map $\psi_0:B\times (G\times S)/\Gamma\to \Lambda$. 
Hence we can define 
a $\beta$-equivariant map $\psi:B^3\times (G\times S)/\Gamma
\to \Lambda^3$ by 
\begin{equation}\psi(a,b,c,x)=
(\psi_0(a,x),\psi_0(b,x), \psi_0(c,x))\, (x\in (G\times S)/\Gamma).\notag
\end{equation}

Let $\Lambda^\prime\subset \Lambda$ be the support of the measure
$(\psi_0)_*(\rho\times \lambda)$.
Then $\Lambda^\prime$ is a closed subset of $\Lambda$. For
every $g\in \Gamma$ and almost every $u\in S$ the element
$\alpha(g,u)\in H$ stabilizes $\Lambda^\prime$, and $\alpha$
is cohomologous to a cocycle with values in the intersection
of $H$ with the stabilizer of $\Lambda^\prime$.
Since the stabilizer of $\Lambda^\prime$ is a closed subgroup 
of ${\rm Iso}(X)$ we may assume without loss of generality that 
$\Lambda^\prime=\Lambda$, i.e. that the support 
of the measure $(\psi_0)_*(\rho\times \lambda)$
equals the limit set of $H$.

For simplicity,
write $L^2(H\times H)$ for the space of functions on $H\times H$
which are square integrable with respect to a Haar measure on $H\times H$,
Let $L^{[2]}((G\times S)/\Gamma,L^2(H\times H))$ be the space of all
measurable maps $(G\times S)/\Gamma\to L^2(H\times H)$ with the additional
property that for each such map $\zeta$ the function $x\to
\Vert\zeta(x)\Vert$ is square integrable on $(G\times S)/\Gamma$
(where $\Vert \,\Vert$ is the $L^2$-norm on $L^2(H\times H)$).
Then $L^{[2]}((G\times S)/\Gamma,L^2(H\times H))$ 
has a natural structure of a
separable Hilbert space, and the group $G$ acts on
$L^{[2]}((G\times S)/\Gamma,L^2(H\times H))$ as a group of isometries. In
other words, $L^{[2]}((G\times S)/\Gamma,L^2(H\times H))$ is a Hilbert
module for $G$.

By Theorem \ref{cocycle},
there is a continuous $L^2(H\times H)$-valued bounded nontrivial cocycle
$\phi:\Lambda^3\to L^2(H\times H)$. The composition of the map
$\psi:B^3\times (G\times S)/\Gamma\to \Lambda^3$ with 
the cocycle $\phi$ can be viewed as a nontrivial 
$\beta$-invariant measurable bounded map $B^3\to
L^{[2]}((G\times S)/\Gamma,L^2(H\times H))$. Since $B$ is a strong
boundary for $G$, this map defines a nontrivial
cohomology class in $H_{cb}^2(G,L^{[2]}((G\times S)/\Gamma,
L^2(H\times H)))$.

Now if $G$ is simple then by the results of Monod
and Shalom \cite{MS04} there is a $\beta$-equivariant map
$(G\times S)/\Gamma\to L^2(H\times H)$. Since the action of $G$ on
$(G\times S)/\Gamma$ is ergodic, by the cocycle reduction lemma of
Zimmer \cite{Z} the cocycle $\beta$ and hence $\alpha$ is
cohomologous to a cocycle into a compact subgroup of $H$.
On the other hand, if $G=G_1\times G_2$ for
semi-simple Lie groups $G_1,G_2$ with finite center and without
compact factors, then the results of Burger and Monod
\cite{BM99,BM02} show that via possibly exchanging $G_1$ and $G_2$
we may assume that there is an equivariant map $(G\times S)/\Gamma
\to L^2(H\times H)$ for the restriction of $\beta$ to $G_1\times
(G\times S)/\Gamma$, viewed as a cocycle for $G_1$. 
By Lemma \ref{mildmixing} the
action of $G_1$ on $(G\times S)/\Gamma$ is ergodic and therefore the
cocycle reduction lemma of Zimmer \cite{Z} shows that the
restriction of $\beta$ to $G_1$ is equivalent to a cocycle into a
compact subgroup of $H$. We now follow the proof of Theorem 1.2 of
\cite{MS04} and find a minimal such compact subgroup $K$ of $H$.
The cocycle $\beta$ and hence $\alpha$ is cohomologous to a cocycle
into the normalizer $N$ of $K$ in $H$.
Moreover, there is a continuous
homomorphism of $G$ onto $N/K$. Since $G$ is connected, the image
of $G$ under this homomorphism is connected as well and hence by
Proposition \ref{firstpart}, $N/K$ is a simple Lie group of rank one.
This completes the proof of the proposition. \end{proof}

We formulate Proposition \ref{cocyclerig} once more in the 
particular case that the measure space $S$ consists
of a single point.

\begin{corollary}\label{homomor} 
Let $G$ be a connected semi-simple Lie
group with finite center, no compact factors and of rank at least
2. Let $\Gamma$ be an irreducible lattice in $G$, let $X$ be a
proper hyperbolic geodesic metric space and let $\rho:\Gamma\to
{\rm Iso}(X)$ be a homomorphism. Let $H<{\rm Iso}(X)$ be the
closure of $\rho(\Gamma)$. If $H$ is non-elementary, then $H$ is a
compact extension of a simple Lie group $L$ of rank one, and if
$\pi:H\to L$ is the canonical projection then $\pi\circ \rho$
extends to a continuous surjective homomorphism $G\to L$.
\end{corollary}

A proper hyperbolic geodesic metric space $X$ is
\emph{of bounded growth} if there is a number $b>0$ such that for
every $x\in X$ and every $R>0$ the closed ball of radius $R$ about
$x$ contains at most $be^{bR}$ disjoint balls of radius 1. This
implies that every elementary subgroup of ${\rm Iso}(X)$ is
amenable (Proposition 8 of \cite{MMS04}).
Thus we can use Corollary \ref{homomor} to show.

\begin{corollary}\label{boundedgr} 
Let $G$ be a connected semi-simple Lie
group with finite center, no compact factors,
and no factors of rank one.
Let $\Gamma$ be an irreducible lattice in $G$, let $X$ be a
proper hyperbolic geodesic metric space of bounded growth
and let $\rho:\Gamma\to
{\rm Iso}(X)$ be a homomorphism. Then the closure
of $\rho(\Gamma)$ in ${\rm Iso}(X)$ is compact.
\end{corollary}
\begin{proof}
Let $G,\Gamma,X$ be as in the corollary and let
$\rho:\Gamma\to {\rm Iso}(X)$ be a homomorphism.
Since $G$ has no factor of rank one by assumption, $G$
does not admit any surjective homomorphism onto a simple
Lie group of rank one. Thus
by Corollary \ref{homomor}, the closure $H$ of $\rho(\Gamma)$ in
${\rm Iso}(X)$ is elementary and hence
amenable since $X$ is of bounded growth by assumption.
By Margulis' normal subgroup
theorem, either the kernel $K$ of $\rho$ has finite index in
$\Gamma$ and $\rho(\Gamma)$ is finite,
or $K$ is finite and the group $\Gamma/K= \rho(\Gamma)$ has
Kazhdan's property $T$. But $\rho(\Gamma)$
is a dense subgroup of $H$ and therefore $H$ has property $T$ if
this is the case for $\rho(\Gamma)$.
In other words, $H$
is an amenable group with
property $T$ and hence $H$ is compact (see \cite{Z} for details).
This shows the corollary.
\end{proof}

\section{Totally disconnected groups of isometries}

In this section we investigate closed subgroups
of the isometry group of a proper hyperbolic geodesic
metric space $X$ which are compact extensions of
totally disconnected groups. We continue to use the assumptions
and notations from Sections 2-5. We
begin with a simple observation on topological properties
of a closed subgroup $G$ of ${\rm Iso}(X)$ 
on its limit set.

\begin{lemma}\label{closedorbit} 
Let $X$ be a proper hyperbolic geodesic
metric space and let 
$G<{\rm Iso}(X)$ be a closed subgroup of ${\rm Iso}(X)$ with
limit set $\Lambda$. Let $\Delta$ be the diagonal in 
$\Lambda\times \Lambda$.
For every hyperbolic element
$g\in G$ with fixed points $a\not=b\in \Lambda$, 
the $G$-orbit of
$(a,b)$ is a closed subset of $\Lambda\times\Lambda-\Delta$.
\end{lemma}
\begin{proof}
Let $X$ be a proper hyperbolic geodesic
metric space and let $G<{\rm Iso}(X)$ be a closed
subgroup with limit set $\Lambda$.
Assume without loss of generality that $G$ is non-elementary.
Let $g\in G$ be a hyperbolic element with fixed points
$a\not=b\in \Lambda$. We have to show that the $G$-orbit
of $(a,b)$ is a closed subset of $\Lambda\times \Lambda-\Delta$.

Thus let $(\xi_i,\eta_i)\subset \Lambda\times \Lambda-\Delta$ 
be a sequence
of pairs of points in this orbit which converges to some
$(a^\prime,b^\prime)\in \Lambda\times \Lambda-\Delta$. Then there
are $u_i\in G$ such that $(u_ia,u_i b)=(\xi_i,\eta_i)$. Let
$\Gamma=\{g^k\mid k\in \mathbb{Z}\}$  be the cyclic subgroup of
$G$ generated by $g$. Let $\gamma$ be a geodesic connecting
$a$ to $b$ and let $x_0=\gamma(0)$. The Hausdorff distance
between $\gamma$ and the orbit $\Gamma x_0$ 
of $x_0$ under the action of
$\Gamma$ is finite, say this Hausdorff distance is not bigger than
some number $k_0>0$.
As a consequence, every
point on the geodesic $u_i\gamma$ is contained in the
$k_0$-neighborhood of the orbit of $u_ix_0$ under the action of 
the infinite cyclic group
$u_i\Gamma u_i^{-1}$. 

Since $(\xi_i,\eta_i)\to (a^\prime,b^\prime)\in
\Lambda\times \Lambda-\Delta$, the geodesics $u_i\gamma$
connecting the points $\xi_i,\eta_i$ intersect a fixed compact subset
$K$ of $X$ which is independent of $i$. 
Thus via composing $u_i$ with
$g^{\ell_i}$ for a suitable number 
$\ell_i\in \mathbb{Z}$ we may assume
that the image of the point $x_0$ under the isometry $u_i$ is
contained in the closed $k_0$-neighborhood $K^\prime$ of the
compact set $K$. Since the subset of ${\rm Iso}(X)$ of all
isometries which map $x_0$ to a fixed compact subset of $X$ is
compact and since $G<{\rm Iso}(X)$ is closed by assumption, after
passing to a subsequence we may assume that the sequence 
$(u_i)\subset G$
converges in $G$ to an element $u\in G$ which maps $(a,b)$ to
$(a^\prime,b^\prime)$. Therefore $(a^\prime,b^\prime)$ is
contained in the $G$-orbit of $(a,b)$. This shows the lemma.
\end{proof}

Following \cite{BM99}, for any locally compact
$\sigma$-compact topological group $G$, every  
element in the kernel of the homomorphism 
$H_{cb}^2(G,\mathbb{R})\to H_c^2(G,\mathbb{R})$ can
be represented by a \emph{continuous quasi-morphism}. Such a
continuous quasi-morphism
is a continuous map $\rho:G\to \mathbb{R}$ such that
\begin{equation}
\sup_{g,h\in G}\vert \rho(g)+\rho(h)-\rho(gh)\vert <\infty. \notag
\end{equation}
The set of all continuous quasi-morphisms for $G$ is naturally
a vector space. 
A continuous quasi-morphism is bounded on 
every \emph{compact} subset of $G$ and is bounded
on each fixed conjugacy class in $G$. The quasi-morphism $\rho$ 
defines a non-trivial element in the kernel of the
natural homomorphisms 
$H_{cb}^2(G,\mathbb{R})\to H_c^2(G,\mathbb{R})$ if
and only if there is no continuous homomorphism
$\chi:G\to \mathbb{R}$ such that $\rho-\chi$ is bounded.

We say that the quasi-morphism $\rho$ \emph{separates}
an element $g\in G$ from a subset $A\subset G$ if
$\lim\inf_{k\to \infty}\frac{1}{k}\vert \rho(g^k)\vert >0$
and $\lim\sup_{k\to \infty}\frac{1}{k}\vert \rho(h^k)\vert =0$
for every $h\in A$.
The following proposition extends earlier results of 
Fujiwara \cite{F98} (see also the papers \cite{BF02,H08} for
generalizations of \cite{F98} in a different direction).

\begin{proposition}\label{quasimorphism}
Let $X$ be a hyperbolic geodesic metric space and let $G<{\rm Iso}(X)$
be a closed non-elementary subgroup with limit
set $\Lambda$. Let $g_1,\dots,g_k\in G$
be any hyperbolic elements and let $a_i,b_i$ be the attracting and
repelling fixed point for the action of $g_i$ on $\Lambda$,
respectively.
If the ordered pairs $(a_i,b_i),(b_j,a_j)\in 
\Lambda\times \Lambda$ $(i,j\leq k)$ are contained in 
pairwise distinct orbits for
the action of $G$ on $\Lambda\times \Lambda$ 
then for every $i$ there is 
a continuous quasi-morphism $\rho_i:G\to \mathbb{R}$
which separates $g_i$ from $\{g_j\mid j\not=i\}$.
\end{proposition}
\begin{proof}
We follow the strategy from the 
proof of Theorem \ref{cocycle}. Let 
$\kappa_0>0$ be a hyperbolicity constant for $X$. Let $\Lambda\subset
\partial X$ be the limit set of a closed subgroup
$G<{\rm Iso}(X)$. Then
for every hyperbolic element $g\in G$ with fixed points
$a,b\in \Lambda$, for every
geodesic $\gamma$ connecting $a$ to $b$ and for every $x_0\in \gamma$,
the orbit $\{g^kx_0\mid k\in \mathbb{Z}\}$ of $x_0$ under the
action of the infinite cyclic subgroup of $G$ generated
by $g$  is contained in the $\kappa_0$-neighborhood of $\gamma$. 

Let $\delta_z$ $(z\in X)$ be a family
of metrics on $\partial X$ which satisfies the properties
(\ref{dist}),(\ref{distequi}) from Section 2.
As in the proof of
Theorem \ref{cocycle}, 
choose $c_0>0$ such that $\delta_z(\xi,\eta)\geq 2c_0$
for all $\xi\not=\eta\in \Lambda$ and every
$z\in X$ whose distance to a geodesic connecting
$\xi$ to $\eta$ is at most $\kappa_0$. 
Let $\phi_0:[0,\infty)\to
[0,1]$ be a smooth function such that
$\phi_0(t)=0$ for all $t\leq c_0$ and $\phi_0(t)=1$ for
all $t\geq 2c_0$ and define a function 
$\phi$ on $V=(\Lambda\times \Lambda-\Delta)\times X\times X$
by $\phi(\xi,\eta,x,y)=
\phi_0(\delta_{x}(\xi,\eta))\phi_0(\delta_{y}(\xi,\eta)).$

Let $\kappa>0$ be sufficiently large that the
set $\{x\in X\mid \delta_x(\xi,\eta)\geq c_0\}$ is contained
in the $\kappa$-tubular neighborhood of any geodesic
connecting $\xi$ to $\eta$. Let
$\psi:[0,\infty)\to [0,1]$ be a smooth function which
satisfies $\psi(t)=1$ for $t\leq \kappa$ and
$\psi(t)=0$ for $t\geq 2\kappa$ and define a function
$\zeta$ on $V$
by $\zeta(\xi,\eta,x,y)=\phi(\xi,\eta,x,y)
\psi(d(x,y))$. The function $\zeta$ is invariant under
the diagonal action of $G$ and hence it projects to
a function on $W=G\backslash V$ whose support we denote by $W_0$.
Let $\hat d$ be the $G$-invariant
metric on $V$ defined as in 
equation (\ref{productmetric2}) in the proof of 
Theorem \ref{cocycle} and let $d_0$ be its projection to
$W=G\backslash V$. 

Define a $G$-invariant
foliation ${\cal F}$ of $V$ by requiring that
the leaf of ${\cal F}$ through $(\xi,\eta,x,y)$ equals
$F(\xi,\eta)=\{(\xi,\eta,x,y)\mid x,y\in X\}$. The foliation
${\cal F}$ then projects to a foliation $\hat {\cal F}$ on $W$.
Let $\iota:V\to V$ be the $G$-equivariant 
involution which maps a point
$(\xi,\eta,x,y)\in V$ to $(\eta,\xi,x,y)$. 
This involution projects to an involution of $W$ which we
denote again by $\iota$. Let $P:V\to W$ be the canonical projection.

Let $g_1,\dots,g_k\in G$ be hyperbolic elements
with attracting and repelling fixed points $a_i,b_i\in \Lambda$,
respectively.
Assume that the $G$-orbits of the ordered pairs
$(a_i,b_i),(a_j,b_j)\in \Lambda\times \Lambda$ 
$(i,j\in \{1,\dots,k\})$ are
pairwise disjoint. Then we have 
$i(PF(a_i,b_i))\cap PF(a_j,b_j)=\emptyset$ for all $i,j$. 
In particular, 
the action of $\iota$ on 
$W_0\cap \bigcup_iPF(a_i,b_i)$ is non-trivial.
Since by Lemma \ref{closedorbit} for every $i\in \{1,\dots,k\}$ 
the $G$-orbit of
$(a_i,b_i)$ is a closed subset of $\Lambda\times \Lambda-\Delta$,
we can find a small ball
$B\subset W$ which is disjoint from its image under $\iota$
and such that
the interior of the preimage $\tilde B$ of $B$ in $V$ meets
the intersection of
$F(a_i,b_i)$ with the support of $\zeta$ and that 
$\tilde B\cap F(a_j,b_j)=\emptyset$ for 
$j\not=i$.

As in the proof of Theorem \ref{cocycle}, 
let ${\cal H}$ be the
vector space of H\"older continuous functions on $W$
supported in $B$. For $f\in {\cal H}$ denote
by $\tilde f$ the lift of $f$ to a $G$-invariant
$\iota$-anti-invariant
function on $V$. For $\xi\not=\eta\in\Lambda$ define
$f_{\xi,\eta}$ 
to be the restriction of $\tilde f\zeta$ to the
leaf $F(\xi,\eta)$, viewed as a function on $X\times X$
(note that this is slightly inconsistent with the notations
in the proof of Theorem \ref{cocycle}).
For a subset $C$ of
$X$ and a number $r>0$ define 
$N(C,r)=\{(x,y)\in X\times X\mid d(x,C)\leq r,
d(y,C)\leq r, d(x,y)\leq 2\kappa\}$.
Then for all $\xi\not=\eta\in \Lambda$ and every geodesic
$\gamma$ connecting $\eta$ to $\xi$ 
the support of $f_{\xi,\eta}$ 
is contained in 
$N(\gamma,\kappa)$.
By Proposition 5.2 of \cite{MS04}, there is
a $G$-invariant Radon measure $\mu$ on $X$ with full
support such that the $\mu$-mass of a ball in
$X$ of radius $4\kappa$ is at most one.
Let $\nu=\mu\times \mu$; it
follows from the discussion in the proof of Theorem \ref{cocycle}
that for every subset $C$ of $X$ of diameter
at most $R$ the $\nu$-mass of the intersection
$N(C,\kappa)\times N(\gamma,\kappa)$ is bounded
from above by $mR$ for a universal constant
$m>0$.

Choose an arbitrary point $\xi\in \Lambda$ and a small
open neighborhood $A$ of $\xi$ in $X\cup \partial X$.
Identifying each leaf of the foliation ${\cal F}$ with 
$X\times X$,
for $f\in {\cal H}$ and an element $g\in G$ with $g\xi\not=\xi$
define
\begin{equation}
\Phi(f)(g)=\int_{F(\xi,g\xi)-A\times A-gA\times gA}
f_{\xi,g\xi}d\nu \notag
\end{equation}
and if $g\xi=\xi$ define $\Phi(f)(g)=0$.
By construction, the intersection of the support
of $f_{\xi,g\xi}\zeta_{\xi,g\xi}$ with $X\times X-A\times A-gA\times gA$
is compact and hence this integral exists.

We claim that $\Phi(f)(g)$ depends continuously on $g$.
Namely, if $(g_i)\subset G$ is any sequence which
converges to some $g\in G$ $(i\to \infty)$ 
then $g_i(X\cup \partial X-A)$ converges to
$g(X\cup\partial X-A)$ in the Hausdorff topology for compact
subsets of $X\cup \partial X$, moreover $g_i\xi\to g\xi$. 
Since $\tilde f\zeta$ is continuous and bounded, 
the functions $f_{\xi,g_i\xi}$ converge
locally in $L^1(X\times X,\nu)$ to $f_{\xi,g\xi}$ and consequently 
$\Phi(f)(g_i)\to \Phi(f)(g)$ $(i\to \infty)$.
In other words, $g\to \Phi(f)(g)$ is continuous.

Next we claim that $\Phi(f)$ is a quasi-morphism
for $G$. More precisely, there is some number $c>0$ only
depending on the H\"older norm of $f$ such that
\begin{equation}
\vert \Phi(f)(g)+\Phi(f)(h)-\Phi(f)(gh)\vert \leq c \notag
\end{equation}
for all $g,h\in G$. To see this simply observe that
by invariance of $\tilde f\zeta$
under the action of $G$ and by anti-invariance of $\tilde f\zeta$
under the involution $\iota$ we have
\begin{align}
\Phi(f)(g)+ \Phi(f)(h)-\Phi(f)(gh) &
=\int_{F(\xi,g\xi)-A\times A-gA\times gA}f_{\xi,g\xi} d\nu  \notag\\
+\int_{F(g\xi,gh\xi)-gA\times gA-ghA\times ghA}f_{g\xi,gh\xi} d\nu &
+\int_{F(ghA,A)-ghA\times ghA-A\times A}f_{gh\xi,\xi} d\nu \notag
\end{align}
and therefore our claim is an immediate
consequence of the estimates in Step 3 and Step 4 of the proof
of Theorem \ref{cocycle}.

For every $j\in \{1,\dots,k\}$  
the infinite cyclic subgroup $\Gamma$ generated by $g_j$ 
acts properly on the
intersection of $F(a_j,b_j)$ with the support of the function 
$\zeta$, with compact fundamental domain $D_j$. By our
choice of $B$, for every $q\in \mathbb{R}$ 
there is a function $f\in {\cal H}$ with
the property that $\int_{D_i} f_{a_i,b_i}d\nu =q$ and
$\int_{D_\ell} f_{a_\ell,b_\ell}d\nu=0$ for $\ell\not=i$.
Now the arguments in the
proof of Theorem 3.1 of \cite{H08} show that 
$\lim_{k\to \infty} \Phi(f)(g_i^k)/k=q$ and 
$\lim_{k\to \infty}\Phi(f)(g_\ell^k)/k=0$ for all
$\ell\not=1$. In other words, $\Phi(f)$
is a continuous quasi-morphism for $G$ which separates
$g_i$ from $\{g_j\mid j\not=i\}$.
\end{proof}

The following
corollary shows the first part of Theorem 2 from the
introduction and extends earlier results of Fujiwara \cite{F98} 
(see also \cite{BF02,H08}).

\begin{corollary}\label{trivialcoeff}
Let $X$ be a proper hyperbolic geodesic
metric space and let $G$ be a closed non-elementary subgroup of
${\rm Iso}(X)$ with limit set $\Lambda\subset \partial X$. If the
action of $G$ on the complement of the diagonal in $\Lambda\times
\Lambda$ is not transitive then the kernel of the natural 
homomorphism $H_{cb}^2(G,\mathbb{R})\to H_c(G,\mathbb{R})$
is infinite dimensional.
\end{corollary}
\begin{proof} Let $G<{\rm Iso}(X)$ be a closed
non-elementary subgroup with limit set $\Lambda$.
Denote by $\Delta$ the diagonal in $\Lambda\times \Lambda$ and
assume that $G$ does not act transitively on $\Lambda\times\Lambda
-\Delta$. 
The set of pairs of fixed points of hyperbolic elements
of $G$ is dense in $\Lambda\times \Lambda-\Delta$,
and the action of $G$ on $\Lambda\times\Lambda-\Delta$ has a 
dense orbit
\cite{GH}. Let $g\in G$ be a hyperbolic element 
and let $(a,b)\in \Lambda\times \Lambda$ be the ordered
pair of fixed points for 
the action of $G$ on $\partial X$.
By Lemma \ref{closedorbit}, the  
$G$-orbit $G(a,b)$ of $(a,b)$ is a closed subset of
$\Lambda\times \Lambda-\Delta$. Since $G$ does not
act transitively on $\Lambda\times\Lambda-\Delta$ by
assumption, the complement of $G(a,b)$ in 
$\Lambda\times\Lambda-\Delta$ contains a pair of fixed points
$(a^\prime,b^\prime)$ 
for a hyperbolic element $h\in G$. The orbit $G(a^\prime,b^\prime)$
of $(a^\prime,b^\prime)$ is distinct from the orbit 
$G(a,b)$ of $(a,b)$.

Let $\gamma,\gamma^\prime$ be
geodesics connecting $b,b^\prime$ to $a,a^\prime$. By  
Lemma \ref{closedorbit} and its proof, for every $m>0$ there is a number 
$R>0$ such that for every subsegment $\eta$ of $\gamma$ of length $R$,
there is no $u\in G$ which maps $\eta$ into the $m$-neighborhood
of $\gamma^\prime$. In other words, the
group $G$ satisfies the assumption in Theorem 1 of \cite{BF02}. As a
consequence of Proposition 2 of \cite{BF02} (whose proof is valid
without the assumption that the space $X$ is a graph or
that the group of isometries is countable), there is a
free subgroup $\Gamma<G$ with two generators 
consisting of hyperbolic elements and with
the following properties.
\begin{enumerate}
\item[i)] For a fixed point $x_0\in X$, the orbit
map $u\in \Gamma\to ux_0\in X$ is a quasi-isometric
embedding.
\item[ii)] There are infinitely many elements $u_i\in \Gamma$
$(i>0)$ with fixed points $a_i,b_i$
such that for all $i>0$ the $G$-orbit of
$(a_i,b_i)\in \Lambda\times\Lambda-\Delta$ is distinct from
the orbit of $(b_j,a_j)$ $(j>0)$ or
$(a_j,b_j)$ $(j\not=i)$.
\end{enumerate}

Choose $\{h_1,\dots,h_n\}\subset \{u_i\mid i>0\}\subset \Gamma$
as in ii) above. By Proposition \ref{quasimorphism},
for every $i$ there is a continuous quasi-morphism
$\rho_i$ for $G$ which separates $h_i$ from $\{h_j\mid
j\not=i\}$. This implies that the dimension of the kernel of 
the natural homomorphism $H_{cb}^2(G,\mathbb{R})\to H_c(G,\mathbb{R})$
is at least $n$. Since $n>0$ was arbitrary we conclude that 
the kernel of the natural 
map $H_{cb}^2(G,\mathbb{R})\to H_c(G,\mathbb{R})$ is 
indeed infinite dimensional. 
\end{proof}

{\bf Remark:} The proof of Corollary \ref{trivialcoeff} also
shows  
the following. If $G<{\rm Iso}(X)$ is a closed subgroup
with limit set $\Lambda$ 
whose action on $\Lambda\times \Lambda-\Delta$  
is not transitive then there is an infinite dimensional
vector space of continuous bounded $G$-invariant functions
$\omega:\Lambda^3\to \mathbb{R}$ which are anti-symmetric under
permutations of the three variables and satisfy the cocycle
equation (\ref{coc}).

\bigskip

The following proposition completes the proof of 
Theorem 2 from the introduction.

\begin{proposition}\label{transitive} 
Let $G<{\rm Iso}(X)$ be a closed
non-elementary subgroup with limit set $\Lambda$. If $G$ acts
transitively on the complement of the diagonal in $\Lambda\times
\Lambda$ then the kernel of the natural homomorphism
$H_{cb}^2(G,\mathbb{R})\to H_c^2(G,\mathbb{R})$ is trivial.
\end{proposition}
\begin{proof} Let $G<{\rm Iso}(X)$ be a closed non-elementary
subgroup which acts transitively on the space $A$ of pairs of
distinct points in $\Lambda$. Since every element 
in the kernel of the natural map $H_{cb}^2(G,\mathbb{R})\to 
H_c^2(G,\mathbb{R})$ can be represented by a continuous
unbounded quasi-morphism it suffices to show that such 
a continuous unbounded quasi-morphism for $G$ does not exist.

Let $a\in \Lambda$ be any point.
For $x,y\in X$ and a geodesic ray
$\gamma:[0,\infty)\to X$ connecting $x$ to $a$ write
$\beta(y,\gamma)=\lim\sup_{t\to \infty}(d(y,\gamma(t))-t)$ and
define the \emph{Busemann function}
\begin{equation}\label{busemann}
\beta_a(y,x)=\sup\{\beta(y,\gamma)\mid \gamma\,\text{is a
geodesic ray connecting}\,x\,\text{to}\,a\}.
\end{equation}
By Lemma 8.1 of \cite{GH}, there is 
a constant $c>0$ with the following properties. Let 
$\gamma:\mathbb{R}\to X$ be any geodesic
with $\gamma(t)\to a$ $(t\to\infty)$. 
Then for every fixed $s\in \mathbb{R}$ and all sufficiently large
$T>0$ we have 
\begin{equation}\label{buseone}
\vert\beta_a(\cdot,\gamma(s))- (d(\cdot,\gamma(T))-T+s)\vert \leq c.
\end{equation}
This implies that
\begin{equation}\label{busetwo}
\vert \beta_a(\cdot,\gamma(s))-\beta_a(\cdot,\gamma(t))+s-t\vert \leq 2c
\end{equation}
for all $s,t\in \mathbb{R}$.
Morever, by Proposition 8.2 of \cite{GH}, for all $x,y\in X$ we have
\begin{equation}\label{buse}
\vert \beta_a(\cdot,y)-\beta_a(\cdot,x)+\beta_a(y,x)\vert \leq c
\end{equation}
and consequently
\begin{equation}\label{busethree}
\vert \beta_a(x,y)-\beta_a(y,x)\vert \leq 2c.
\end{equation}

Define the \emph{horosphere} at $a$
through $x$ to be the set 
\begin{equation}\label{horosphere}
H_a(x)=\beta_a(\cdot,x)^{-1}[-4c,4c].
\end{equation}
By inequality (\ref{busethree}), if $z\in X$ is any
point with $\vert \beta_a(x,z)\vert \leq 2c$
then $\vert\beta_a(z,x)\vert \leq 4c$ and hence
$z\in H_a(x)$.

We claim that there is a universal constant
$c_0>0$ with the following property.
Let $\gamma:\mathbb{R}\to X$ be a
biinfinite geodesic with $\gamma(t)\to a$
$(t\to \infty)$. Then for every point $y\in X$
and every $t\in \mathbb{R}$ we have
\begin{equation}\label{horodistance}
d(\gamma(t),H_a(y))\leq
\vert\beta_a(y,\gamma(t))\vert.
\end{equation}

To see this let $p=\beta_a(y,\gamma(t))$. 
Then the estimate (\ref{busetwo}) shows that
$\vert\beta_a(y,\gamma(t+p))\vert \leq 2c$.
Thus by the inequality (\ref{busethree})
we have $\vert \beta_a(\gamma(t+p),y)\vert \leq 4c$ and
hence $\gamma(t+p)\in H_a(y)$. This shows the claim

Since $G$ is non-elementary by assumption, $G$ contains a
hyperbolic element $g\in G$. Let $a\in \Lambda$ be the attracting
fixed point of $g$ and let $b\in \Lambda-\{a\}$ be the repelling
fixed point. Then $g$ preserves the set of geodesics connecting
$b$ to $a$.
Let $W(a,b)\subset X$
be the closed non-empty subset of all points in $X$ which lie on
a geodesic connecting $b$ to $a$. The set $W(a,b)$ is contained
in a tubular neighborhood of fixed radius $\kappa_0>0$
about a fixed geodesic connecting $b$ to $a$. 
The isometry $g\in G$ is hyperbolic with fixed points $a,b\in \Lambda$
and therefore it preserves $W(a,b)$.
If we denote by $\Gamma$ the infinite cyclic
subgroup of $G$ generated by $g$ then $W(a,b)/\Gamma$ is
compact.
As a consequence of the estimate (\ref{buseone}),
there is a number $\nu>0$ and for every
$x\in W(a,b)$
and every $t\in \mathbb{R}$ there is a number $k(t)\in
\mathbb{Z}$ with $\vert \beta_a(g^{k(t)}(x),x)-t\vert <\nu$. It
follows from this and (\ref{buse})
that  for every $y\in X$ there is some $k=k(y)\in
\mathbb{Z}$ such that $\vert \beta_a(y,g^k(x))\vert \leq \nu+c$. 
Since $g^k(x)$ is contained in a biinfinite geodesic converging
to $a$, we obtain from the estimate (\ref{horodistance}) above 
that the distance between $g^k(x)$ and the horosphere
$H_a(y)$ is at most $\delta_0=\nu+c$.

Since $G$ acts transitively on the complement of the
diagonal in $\Lambda\times \Lambda$, 
the stabilizer $G_a<G$ of the point $a\in \Lambda$
acts transitively on $\Lambda-\{a\}$. Thus there is for every
$\zeta\in \Lambda-\{a\}$ an element $h_\zeta\in G_a$ with
$h_\zeta(b)=\zeta$. Let again $x\in W(a,b)$.
Since $h_\zeta\in G_a$ we
have $h_\zeta^{-1}(H_a(x)) =H_a(h_\zeta^{-1}(x))$. By our above
consideration, there is a number $\ell\in \mathbb{Z}$ such that
the distance between $g^{\ell}(x)$ and $H_a(h_\zeta^{-1}(x))$
is at most $\delta_0$ and therefore the
distance between $h_\zeta\circ g^\ell(x)$ and
$H_a(x)$ is at most $\delta_0$. Hence via
replacing $h_\zeta$ by $h_\zeta\circ g^\ell\in G_a$ we may assume
that the distance between $h_\zeta(x)$ and $H_a(x)$
is at most $\delta_0$.

For $x\in X$
denote by $N_{a,x}\subset G_a$ the set of all elements $h\in G_a$ with
the property that the distance between $H_a(x)$ and
$hx$ is at most $\delta_0$. The above consideration
shows that for every $x\in W(a,b)$
and every $\zeta\in \Lambda-\{a\}$ there is
some $h_\zeta\in N_{a,x}$ which maps $b$ to $\zeta$.

For every $h\in N_{a,x}$ the sequence
$(g^{-\ell}\circ h\circ g^{\ell})_{\ell>0}$ is contained in $G_a$.
We claim that there is a number $\delta_1>0$ not
depending on $h$ such that 
$(g^{-\ell}\circ h\circ g^{\ell})(x)$ is contained in 
the $\delta_1$-neighborhood of $H_a(x)$ for every $\ell >0$.
Namely, inequality (\ref{buseone}) together
with the triangle inequaltity shows that
\begin{equation}\label{busefour}
\vert \beta_a(y,x)-\beta_a(z,x)\vert \leq d(y,z)+2c \notag
\end{equation}
for all $z,y\in X$ and hence 
we conclude that 
$\vert \beta_a(h(x),x)\vert \leq \delta_0+6c$ for
every $h\in N_{a,x}$.  As a consequence
of this and 
inequality (\ref{buse}), we obtain that  
\begin{equation}\label{busecomp}
\vert \beta_a(\cdot,x)-\beta_a(\cdot, h(x))\vert \leq \delta_0+7c.
\notag
\end{equation}

Let again $\gamma:\mathbb{R}\to X$ be a geodesic connecting
$b$ to $a$ with $\gamma(0)=x$. Then 
$h$ maps the geodesic $\gamma$ to a
geodesic $h(\gamma)$ connecting $h(b)$ to $a$.
By inequality (\ref{busetwo}), we have
$\vert \beta_a(\cdot,h(\gamma(t)))-\beta_a(\cdot,h(x))+t\vert
\leq 2c$ for all $t\in \mathbb{R}$ 
and the same estimate for the geodesic $\gamma$ then implies that
\begin{equation}
\vert \beta_a(\cdot,\gamma(t))-\beta_a(\cdot,h(\gamma(t)))\vert \leq
\delta_0+12c \text{ for all } t\in \mathbb{R}. \notag
\end{equation}
Together with the estimate (\ref{horodistance}) we deduce that
for every $t>0$
the distance between $h(\gamma(t))$ and
$H_a(\gamma(t))$ is bounded from above by $\delta_0+12c$.
Since the infinite cyclic group 
$\Gamma$ preserves the set of geodesics connecting
$b$ to $a$ we conclude that there is indeed a universal constant
$\delta_1>0$ such that for every $\ell\geq 0$ the distance
between $h(g^\ell x)$ and $H_a(g^\ell(x))$ is at most
$\delta_1$. This shows the above claim.

Now for every $h\in N_{a,x}$ the sequence
$(g^{-\ell}\circ h\circ g^\ell(b))_{\ell>0}\subset \Lambda$ converges
as $\ell\to \infty$ to $b$. This means that there is a number
$\delta_2>\delta_1$ such that
for sufficiently large
$\ell$ the element $g^{-\ell}\circ h\circ g^{\ell}\in G_a$ maps
the point $x$ into the closed $\delta_2$-neighborhood $B_x$ of $x$.
The group
$\Gamma$ acts on $W(a,b)$ cocompactly and hence if
$C\subset W(a,b)$ is a compact
fundamental domain for this action, then
$B=\cup_{x\in C}B_x$ is compact. Moreover, for every $x\in W(a,b)$,
every element $h\in N_{a,x}$
is conjugate in
$G$ to an element in the compact subset $K=\{u\in G\mid
uB\cap B\not=\emptyset\}$ of $G$.
As a consequence, the restriction to
$N_{a,b}=\cup_{x\in  W(a,b)}N_{a,x}$ of any continuous
quasi-morphism $q$ on $G$ is uniformly bounded.
By our
assumption on $G$ the sets $N_{a,b}$ $((a,b)\in A)$ are pairwise
conjugate in $G$ and hence $q$ is uniformly bounded on
$\cup_{(a,b)\in A}N_{a,b}=N$.

Next we show that the restriction of a quasi-morphism $q$ to the
subgroup $G_{a,b}$ of $G_a$ which stabilizes both points $a,b\in
\Lambda$ is bounded. For this consider an arbitrary element $u\in
G_{a,b}$. We may assume that $u\not\in N_{a,b}$.
Let $x\in C\subset W(a,b)$ where as before, $C$ is a compact
fundamental domain for the action of the infinite
cyclic group $\Gamma$ on $W(a,b)$.

By assumption, $G$ acts transitively on the complement of the
diagonal in $\Lambda\times \Lambda$ and hence
there is an element $h\in G$ such that $h(a)=b$ and $h(b)=a$.
Then $hW(a,b)=W(a,b)$ and hence via composition of $h$ with
and element of $\Gamma$ we may assume that $hx\in C$.
Thus if $\delta_3>0$ is the diameter of $C$ then
\begin{equation}
\vert \beta_b(hux,x)-\beta_a(ux,x)\vert\leq
\vert\beta_b(hux,hx)-\beta_a(ux,x)\vert +c
+d(hx,x)\leq \delta_3+c\notag
\end{equation}
by the estimate (\ref{buse}).

On the other
hand, from another application of the estimate (\ref{buseone}) we
obtain the existence of a constant $\delta_4>0$ such that
$\vert\beta_b(z,y)-\beta_a(y,z)\vert \leq \delta_4$ for all
$y,z\in W(a,b)$. But this just means that the distance between
$hux$ and $u^{-1}x$ is uniformly bounded and hence the
distance between $uhux$ and $x$ is uniformly bounded as
well. As a consequence, the element $uhu$ is contained in a
fixed compact subset of $G$. As before, we conclude from this that
$\vert q(u)\vert $ is bounded from above by a universal constant
not depending on $u$ and hence the restriction of $q$ to $G_{a,b}$
is uniformly bounded.

Now let $h\in G$ be arbitrary and assume that $ha=x,hb=y$. We
showed above that there are $h_y\in N\cap G_x,h_x\in N\cap G_b$
with $h_y(y)=b$ and $h_x(x)=a$; then $h^\prime=h_xh_yh\in G_{a,b}$
and $\vert q(h^\prime)-q(h)\vert $ is uniformly bounded. As a
consequence, $q$ is bounded and hence it defines the trivial
bounded cohomology class. This completes the proof of the
proposition. 
\end{proof}

{\bf Examples:} 

1) Let $T$ be a regular $k$-valent tree for some $k\geq 3$.
Then $T$ is a proper hyperbolic geodesic metric space, and its 
isometry group $G$ is totally disconnected. If $\partial T$
denotes the Gromov boundary of $T$ then 
the group $G$ acts 
transitively on the space $Y$ of \emph{triples} of pairwise
distinct points in $\partial T$. Moreover, there is a 
$G$-invariant measure class $\lambda$ on $\partial T$ with the
property that $(\partial T,\lambda)$ is a strong
boundary for $G$ (see \cite{A94,K03}). As a consequence,
every bounded cohomology class $\omega\in H_{cb}^2(G,\mathbb{R})$
can be represented by a $G$-invariant $\lambda^3$-measurable
bounded function $\omega:Y\to \mathbb{R}$ \cite{BM02,M}
which is anti-symmetric under permutations of the three
variables and which satisfies the cocycle condition
(\ref{coc}). Since the action of $G$ on $Y$ is transitive,
such a function has to vanish. In other words, 
$H_{cb}^2(G,\mathbb{R})=\{0\}$. 

2) Let $G$ be a simple rank-one Lie group of non-compact type.
Then $G$ is the isometry group of a negatively curved
symmetric space $X$. The limit set of $G$ is the full
Gromov boundary $\partial X$ of $X$, moreover the
action of $G$ on the complement of the diagonal
in $\partial X\times \partial X$ is transitive.
By Proposition \ref{transitive}
and well known
results on the usual continuous cohomology of $G$, 
the second bounded 
cohomology group $H_{cb}(G,\mathbb{R})$ 
is trivial if $G\not=SU(n,1)$ and
equals $\mathbb{R}$ for $G=SU(n,1)$ for some $n\geq 2$.

\section{Proper hyperbolic spaces of bounded growth}

In this section we investigate more restrictively
proper hyperbolic geodesic metric spaces of bounded growth. 
This means that there is a number $b>1$ such that for every $R>1$,
every metric ball of radius $R$ contains at most $be^{bR}$
disjoint metric balls of radius $1$.
The following proposition is Theorem 4 from the introduction.

\begin{proposition}\label{boundedgrowth} 
Let $\Gamma$ be a finitely generated
group which admits a proper isometric action on a proper hyperbolic
geodesic metric space of bounded growth. If 
$H_b^2(\Gamma,\mathbb{R})$ or
$H_{cb}^2(\Gamma,\ell^2(\Gamma))$ is finite dimensional 
then $\Gamma$ is virtually
nilpotent.
\end{proposition}
\begin{proof}
Let $\Gamma$ be a finitely generated group which
admits a proper isometric action on a proper hyperbolic geodesic
metric space $X$ of bounded growth. Assume
that $\Gamma$ is infinite and that
$H_b^2(\Gamma,\mathbb{R})$ or $H_b^2(\Gamma,\ell^2(\Gamma))$ is
finite dimensional. By the results of Fujiwara \cite{F98} 
(for real coefficients) and by
\cite{H08} (for coefficients $\ell^2(\Gamma)$),
the subgroup $\Gamma$ of ${\rm Iso}(X)$ is elementary.
We have to show that $\Gamma$ is virtually nilpotent. Since
$\Gamma$ is infinite and acts properly on $X$ by assumption, the
limit set of $\Gamma$ is nontrivial and hence it consists of one
or two points. Assume first that the limit set consists of a
single point $a\in
\partial X$. Then $a$ is a fixed point for the action of $\Gamma$
on $\partial X$.

We recall some notations from the proof
of Proposition \ref{transitive}. Namely,
since $X$ is locally compact, for every $x\in X$  there is a
geodesic ray $\gamma:[0,\infty)\to X$ connecting $x$ to $a$.
For $x\in X$ and $a\in \partial X$ let $y\to \beta_a(y,x)$
be the Busemann function determined by $a$ and $x$ as 
defined in (\ref{busemann}) in the proof of Proposition \ref{transitive}.
Recall also from (\ref{buse}) that there is a number 
$c>0$ such that
$\vert\beta_a(\cdot,y)-\beta_a(\cdot,x)+
\beta_a(y,x)\vert \leq c$ for all $x,y\in X$.
For $x\in X$ let $H_a(x)$ be the horosphere through $x$
and $a$ defined in (\ref{horosphere}).
By the estimate (\ref{busefour}), for all $x,y\in X$ and 
all $r\in \mathbb{R}$ the distance between $\beta_a(\cdot,x)^{-1}(r)$
and $y$ is at least $\vert \beta_a(y,x)-r\vert -2c$.
Together with the estimate (\ref{busethree}), this implies that
the distance in $X$ between any two horospheres
$H_a(x),H_a(y)$ for $x,y\in X$ is at least $\vert
\beta_a(y,x)\vert-12c$.

We claim that there is a number $\kappa >0$ only depending on the
hyperbolicity constant of $X$ such that for every $x\in X$ the
group $\Gamma$ maps the horosphere $H_a(x)$ into
$\beta_a(\cdot,x)^{-1}[-\kappa,\kappa].$ To see this let $g\in
\Gamma, x\in X$ be arbitrary. Then $gH_a(x)=H_a(gx)$ and therefore 
if for some $r>4c$ the image of $H_a(x)$
under $g$ is not contained in $\beta_a(\cdot,x)^{-1}[-r,r]$ then
the distance between $H_a(x)$ and $gH_a(x)$ is at least $r-8c$.

Thus if for \emph{every} $x\in X$ the image of $H_a(x)$ under $g$
is \emph{not} contained in $\beta_a(\cdot,x)^{-1}[-r,r]$ then the
minimal displacement $\inf\{d(y,gy)\mid y\in X\}$ of $g$ is at
least $r-8c$. However, by Proposition 8.24 in \cite{GH} there is a
universal constant $r_0>0$ such that every isometry of $X$ whose
minimal displacement is at least $r_0-8c$ is hyperbolic and hence it
generates an infinite cyclic group of isometries of $X$ whose
limit set consists of two points. Since by assumption 
the limit set of $\Gamma$ consists of a unique point,
the group
$\Gamma$ can not contain such an element. As a consequence, for
every $g\in \Gamma$ there is a point $x(g)\in X$ such that
\begin{equation}
gH_a(x(g))\subset \beta_a(\cdot,x(g))^{-1}[-r_0,r_0].\notag
\end{equation}
In
particular, we have $\vert \beta_a(g(x(g)),x(g))\vert \leq r_0$.

From inequality (\ref{buse}) above with $x=x(g),y=g(x(g))$ we
conclude that \begin{equation}
\vert \beta_a(\cdot,x(g))-\beta_a(\cdot,g(x(g)))\vert\leq r_0+c.
\notag \end{equation} 
On the
other hand, inequality (\ref{buse}) applied to an arbitrary point
$y\in X$ and to $x=x(g)$ shows that $\vert
\beta_a(\cdot,y)-\beta_a(\cdot,x(g))+\beta_a(y,x(g))\vert \leq c$
and similarly $\vert \beta_a(\cdot, gy)-\beta_a(\cdot,
g(x(g)))+\beta_a(gy,g(x(g)))\vert \leq c$. Since
$g$ fixes the point $a$ we have
$\beta_a(y,x(g))=\beta_a(gy,g(x(g)))$ and therefore $\vert
\beta_a(\cdot,y)-\beta_a(\cdot,gy)\vert \leq r_0+3c$ for
every $y\in X$. Since $\beta_a(gz,gy)=\beta_a(z,y)$ for all
$g\in \Gamma$, all $y,z\in X$  
this shows
the above claim.

Let $x\in X$ be an arbitrary point and let $\kappa >0$ be such
that for every $g\in \Gamma$ the horosphere $gH_a(x)$ is contained
in $\beta_a(\cdot,x)^{-1}[-\kappa,\kappa]$. Choose a finite
symmetric generating set $g_1,\dots,g_{2k}$ for $\Gamma$; such a
set exists since $\Gamma$ is finitely generated by assumption. Let
$q=\max\{d(x,g_ix)\mid i=1,\dots, 2k\}$. Since for every $t\in
\mathbb{R}$ the distance between $\beta_a(\cdot,x)^{-1}(t)$ and
$H_a(x)$ is at least $\vert t\vert -6c$, each of the points $g_ix$
can be connected to $x$ by a curve of length at most $q$ which is
contained in $\beta_a(\cdot,x)^{-1}[-q-6c,q+6c]$. Similarly, by the
choice of $\kappa$, for every $g\in \Gamma$ the points $gx,gg_ix$
can be connected by a curve of length at most $q$ contained in
\[V=\beta_a(\cdot,x)^{-1}[-\kappa-q-6c,\kappa+q+6c].\] 

For $g\in \Gamma$
let $\vert g\vert$ be the \emph{word norm} of $g$ with respect to
the generating set $g_1,\dots,g_{2k}$, i.e. $\vert g\vert$ is the
minimal length of a word in $g_1,\dots,g_{2k}$ representing $g$.
Then for every $g\in \Gamma$ the point $gx$ can be connected to
$x$ by a curve which is contained in $V$ and whose length is at
most $q\vert g\vert$.

For $y,z\in V$ define $\delta(y,z)\in [0,\infty]$ to be the
infimum of the lengths of any curve in $X$ which connects $y$ to
$z$ and is contained in $V$. Then the restriction of $\delta$ to a
path-connected component $W$ of $V$ containing the $\Gamma$-orbit
of $x$ is a distance which is not smaller than the restriction of
$d$. For $y\in W$ let $B_W(y,R)\subset W$ be the $\delta$-ball of
radius $R$ about $y$. Let $\gamma:[0,\infty) \to X$ be a geodesic
ray connecting $x$ to $a$. We claim that there is a number
$\chi\geq 1$ such that for every $R>0$ the $\delta$-ball
$B_W(x,e^R)$ of radius $e^R$ about $x$ is contained in the ball in
$X$ of radius $\chi R +\chi$ about $\gamma(\chi R)$.

To see this we argue as in \cite{A94}. Recall from hyperbolicity
that there is a number $\kappa_0 \geq 1$ such that for every
geodesic triangle $\triangle\subset X$ with sides $a,b,c$, the
side $c$ is contained in the $\kappa_0$-neighborhood of $a\cup b$.
Let $y\in W$ be such that $\delta(x,y)\leq 2^m$ for some $m\geq
0$. Then there are $2^m+1$ points $x_0=x,\dots,x_{2^m}=y\in W
\subset V$ such that $d(x_i,x_{i+1})\leq 1$ for all $i$. 
Write $\tilde q=\kappa+q+1+12c$. Then for
$i\leq 2^m$ a geodesic $\gamma_{i,1}$ in $X$ connecting $x_{i-1}$
to $x_{i}$ is contained in
$\beta_a(\cdot,x)^{-1}[-\tilde q,\infty)$. 

For each $i\leq
2^{m-1}$ connect the points $x_{2i-2}$ and $x_{2i}$ by a geodesic
$\gamma_{i,2}$. Then $\gamma_{i,2}$ is contained in the
$\kappa_0$-neighborhood of $\gamma_{2i-1,1}\cup \gamma_{2i,1}$ and
hence 
\begin{equation}\gamma_{i,2}\subset \beta_a(\cdot,x)^{-1}
[-\tilde q-\kappa_0-12c,\infty).
\notag
\end{equation}
Write $\alpha=\kappa_0+12c$. By
induction, for $j\leq m$ choose a geodesic $\gamma_{i,2^j}$
connecting the points $x_{(i-1)2^j}$ and $x_{i2^j}$
$(i=1,\dots,2^{m-j})$. A successive application of our argument
implies that $\gamma_{i,2^j}$ is contained in
$\beta_a(\cdot,x)^{-1}[-j\alpha-\tilde q,\infty)$. In particular,
a geodesic $\zeta$ connecting $x$ to $y$ is contained in
$\beta_a(\cdot,x)^{-1}[-m\alpha-\tilde q,\infty)$. On the other
hand, $\zeta$ is a side of an (ideal) geodesic triangle with
vertices $x,y,a$ and the given geodesic ray $\gamma$ as a second
side. Since for $t\geq 0$ we have $\vert
\beta_a(\gamma(t),x)-t\vert\leq c$ by inequality (\ref{buseone}),
we conclude that for any geodesic ray
$\eta$ connecting $y$ to $a$ the distance between $\eta(\alpha m)$
and $\gamma(\alpha m)$ is bounded by a constant $\tilde c>0$ only
depending on the hyperbolicity constant of $X$. As a consequence,
every point $y\in W$ with $\delta(x,y)\leq 2^m$ is contained in
the ball of radius $\alpha m +\tilde c$ about $\gamma(\alpha m)$.
This is just the statement of our claim with
$\chi=\max\{\alpha\log 2, \tilde c\}$.

From this observation we conclude that the group $\Gamma$ has
\emph{polynomial growth}. By definition, this means that there is
a number $p>0$ such that the number of elements $g\in \Gamma$ of
word norm at most $\ell$ is not bigger than $p\ell^{p}$. Namely,
we observed above that the image of $x$ under an element $g\in
\Gamma$ of word norm at most $\ell$ is contained in the ball
$B_W(x,q\ell)$ where $q>0$ is as above. Since
$\Gamma$ acts properly and isometrically on $X$ 
by assumption, there is a number
$j>0$ such that there are at most $j$ elements $g\in \Gamma$ with
$gx\in B(x,4)$. Then for every $z\in X$ the ball
$B(z,2)$ contains at most $j$ points from the orbit $\Gamma x$ of
$x$ counted with multiplicity. Thus the number of elements of
$\Gamma$ of word norm at most $\ell$ is not bigger than $j$ times
the maximal number of disjoint balls of radius $1$ contained in
$B_W(x,q\ell+1)$ where $q>0$ is as above. Namely, if this 
number equals $k>0$ then there are $k$ balls of radius
$2$ which cover $B_W(x,q\ell+1)$.

On the other hand, by our above observation the ball
$B_W(x,q\ell+1)$ is contained in a ball of radius
$\chi(\log(q\ell))+\chi$ in $X$ for a universal number $\chi>0$.
By assumption, a ball in $X$ of radius $R>1$ contains at most
$be^{bR}$ disjoint balls of radius $1$, where $b>0$ is a universal
constant. Together this shows that $\Gamma$ is indeed of
polynomial growth. By a well known result of Gromov, groups of
polynomial growth are virtually nilpotent. This shows the
proposition in the case that the limit set of $\Gamma$ consists of a
single point.

If the limit set $\Lambda$ of $\Gamma$ consists of two distinct
points $a\not=b\in\partial X$ then by invariance of $\Lambda$
under the action of $\Gamma$, every element of $\Gamma$ maps a
geodesic in $X$ connecting these two points to a geodesic with the
same properties. Now the Hausdorff distance between any two such
geodesics is bounded from above by a universal constant. Since the
action of $\Gamma$ on $X$ is moreover proper, this implies
immediately that the group $\Gamma$ is of polynomial growth and
hence virtually nilpotent. \end{proof}

{\bf Example:} Proposition \ref{boundedgrowth} 
is easily seen to be false for
proper hyperbolic geodesic metric spaces which are not of bounded
growth. Namely, let $(M,g)$ be a symmetric space of non-compact
type and higher rank and consider the space $N=M\times \mathbb{R}$
with the warped product metric $e^{2t}g\times dt$. Then $N$ is a
complete simply connected Riemannian manifold whose curvature is
bounded from above by a negative constant and hence $N$ is a
proper hyperbolic geodesic metric space. However, the isometry
group of $N$ contains an elementary subgroup which is a
semi-simple Lie group of non-compact type and higher rank. Note
that the curvature of $N$ is not bounded from below and $N$ is not
of bounded growth.

\bigskip

\noindent
MATHEMATISCHES INSTITUT DER UNIVERSIT\"AT BONN\\
BERINGSTRA\SS{}E 1\\
D-53115 BONN\\

\smallskip

\noindent
e-mail: ursula@math.uni-bonn.de


\begin{thebibliography}{BM99}


\bibitem{A94} S.~Adams, {\em Boundary amenability for word
hyperbolic groups and an application to smooth dynamics of
simple groups}, Topology 33 (1994), 765--783.

\bibitem{A96} S.~Adams, {\em Reduction of cocycles
with hyperbolic targets,} Erg. Th. \& Dyn. Sys. 16 (1996),
1111--1145.


\bibitem{BF02} M. Bestvina, K. Fujiwara, {\it Bounded cohomology
of subgroups of mapping class groups}, Geometry \& Topology 6
(2002), 69--89.



\bibitem{BH} M.~Bridson, A.~Haefliger, {\sl Metric
spaces of non-positive curvature}, Springer, Berlin Heidelberg
1999.






\bibitem{BM99} M.~Burger, N.~Monod, {\it Bounded cohomology
of lattices in higher rank Lie groups}, J. Eur. Math. Soc. 1
(1999), 199--235.

\bibitem{BM02} M.~Burger, N.~Monod, {\it Continuous bounded
cohomology and applications to rigidity theory},
Geom. Funct. Anal. 12 (2002), 219--280.


\bibitem{F98} K. Fujiwara, {\it The second bounded cohomology
of a group acting on a Gromov hyperbolic space}, Proc. London
Math. Soc. 76 (1998), 70--94.




\bibitem{GH} E.~Ghys, P.~de la Harpe, {\sl Sur les
groupes hyperboliques d'apr\`{e}s Mikhael Gromov,}
Birkh\"auser, Boston 1990.






\bibitem{H08} U.~Hamenst\"adt, {\em Bounded cohomology
and isometry groups of hyperbolic spaces}, arXiv:math.GR/0507097,
to appear in J. Eur. Math. Soc.




\bibitem{K03} V.~Kaimanovich, {\em Double ergodicity
of the Poisson boundary and applications to bounded
cohomology}, Geom. Funct. Anal. 13 (2003), 852--861.






\bibitem{MMS04} I.~Mineyev, N.~Monod, Y.~Shalom,
{\it Ideal bicombings for hyperbolic groups and applications},
Topology 43 (2004), 1319--1344.

\bibitem{M} N.~Monod, {\sl Continuous bounded cohomology
of locally compact groups}, Lecture Notes in Math. 1758,
Springer 2001.

\bibitem{MS04} N.~Monod, Y.~Shalom, {\em Cocycle superrigidity
and bounded cohomology for negatively curved spaces}, J. Diff. Geom.
67 (2004), 395--456.

\bibitem{MS05} N.~Monod, Y.~Shalom,
{\em Orbit equivalence rigidity and bounded
cohomology}, Ann. Math. 164 (2006), 825--878..

\bibitem{SW82} K.~Schmidt, P.~Walters, 
{\em Mildly mixing actions of locally compact groups},
Proc. London Math. Soc. 45 (1982), 506--518.

\bibitem{S00} Y.~Shalom, {\em Rigidity of commensurators
and irreducible lattices}, Invent. Math. 141 (2000), 1--54.

\bibitem{Z} R.~Zimmer, {\sl Ergodic theory and semisimple groups},
Birkh\"auser, Boston 1984.


\end{thebibliography}
\end{document}